\documentclass[11pt,a4paper]{article}
\usepackage{authblk}
\usepackage[english]{babel}
\usepackage[T1]{fontenc}
\usepackage[utf8]{inputenc}
\usepackage{amsmath,amsthm,bm,mathrsfs}
\usepackage{mathtools}
\usepackage{bbold}

\usepackage{palatino}
\usepackage{eurosym}
\usepackage{enumerate}
\usepackage[symbol]{footmisc}

\usepackage[toc,page]{appendix}
\usepackage[inline]{enumitem}
\usepackage[%  
    colorlinks=true,
    pdfborder={0 0 0},
    linkcolor=red
]{hyperref}
\usepackage{hyperref}
\hypersetup{colorlinks,linkcolor={blue},citecolor={blue},urlcolor={red}}  
\usepackage[top=2.5cm, left=2.5cm, right=2.5cm, bottom=2.5cm]{geometry}
\usepackage[square,numbers]{natbib}
\numberwithin{equation}{section}
\newcommand {\argmin}   {\text{{\rm argmin}}}
\usepackage{lineno}
\usepackage{comment}
\modulolinenumbers[5]
\newcommand{\beginsupplement}{%
        \setcounter{table}{0}
        \renewcommand{\thetable}{S\arabic{table}}%
        \setcounter{figure}{0}
        \renewcommand{\thefigure}{S\arabic{figure}}%
     }

\numberwithin{equation}{section}

\theoremstyle{definition}
\newtheorem{assumption}{Assumption}[section]

\title{Discrete and continuum models for the coevolutionary dynamics between CD8+ cytotoxic T lymphocytes and tumour cells}
\date{\vspace{-6ex}}
\vspace{3ex}
\author[1]{Luis Almeida}
\author[2,3]{Chloe Audebert}                             
\author[1]{Emma Leschiera\thanks{Corresponding author \\
\textit{Email addresses}:
 \texttt{almeida@ljll.math.upmc.fr} (Luis Almeida), \texttt{chloe.audebert@sorbonne-universite.fr} (Chloe Audebert),\texttt{leschiera@ljll.math.upmc.fr} (Emma Leschiera), \texttt{tommaso.lorenzi@polito.it} (Tommaso Lorenzi).\\
E.L. has received funding from the European Research Council (ERC) under the European Union’s Horizon2020 research and innovation programme (grant agreement No 740623). \\
T.L. gratefully acknowledges support of the MIUR grant ``Dipartimenti di Eccellenza 2018-2022''.}}
\author[4]{Tommaso Lorenzi$^*$}

\affil[1]{\small\textit{Sorbonne Universit\'e, CNRS, Universit\'e de Paris, Inria, Laboratoire Jacques-Louis Lions UMR 7598, 75005 Paris, France.}}

\affil[2]{\textit{Sorbonne Universit\'e, CNRS, Universit\'e de Paris, Laboratoire Jacques-Louis Lions UMR 7598, 75005 Paris, France.}}

\affil[3]{\textit{Sorbonne Universit\'e, CNRS, Institut de biologie Paris-Seine (IBPS), Laboratoire de Biologie Computationnelle et Quantitative UMR 7238, 75005 Paris, France.}}
\affil[4]{\textit{Department of Mathematical Sciences ``G. L. Lagrange'', Dipartimento di Eccellenza 2018-2022, Politecnico di Torino, 10129 Torino, Italy.}}
\begin{document}
%\begin{frontmatter}
\maketitle
\begin{abstract}
We present an individual-based model for the coevolutionary dynamics between CD8+ cytotoxic T lymphocytes (CTLs) and tumour cells. In this model, every cell is viewed as an individual agent whose phenotypic state is modelled by a discrete variable. For tumour cells this variable represents a parameterisation of the antigen expression profiles, while for CTLs it represents a parameterisation of the target antigens of T-cell receptors (TCRs). We formally derive the deterministic continuum limit of this individual-based model, which comprises a non-local partial differential equation for the phenotype distribution of tumour cells coupled with an integro-differential equation for the phenotype distribution of CTLs. The biologically relevant homogeneous steady-state solutions of the continuum model equations are found. The linear-stability analysis of these steady-state solutions is then carried out in order to identify possible conditions on the model parameters that may lead to different outcomes of immune competition and to the emergence of patterns of phenotypic coevolution between tumour cells and CTLs. We report on computational results of the individual-based model, and show that there is a good agreement between them and analytical and numerical results of the continuum model. These results shed light on the way in which different parameters affect the coevolutionary dynamics between tumour cells and CTLs. Moreover, they support the idea that TCR-tumour antigen binding affinity may be a good intervention target for immunotherapy and offer a theoretical basis for the development of anti-cancer therapy aiming at engineering TCRs so as to shape their affinity for cancer targets.
\end{abstract}


\vspace*{15pt}
\section{Introduction}
\label{sec:intro}
\paragraph{Essentials of the underlying biological problem} CD8+ cytotoxic T lymphocytes (CTLs) play a key role in the immune response against cancer. CTLs carry a specific receptor on their surface, the T-cell receptor (TCR), which can recognise and bind to non-self antigens expressed by tumour cells~\citep{doi:10.1634/theoncologist.2015-0282}. Each TCR recognises and binds specifically to a certain antigen (i.e. the cognate antigen)~\citep{coulie2014tumour}, and possibly other antigens within a certain affinity range~\citep{mason1998, wooldridge2012single}. This enables CTLs to exert an antigen-specific cytotoxic activity against tumour cells, whose efficacy may depend on the affinity range of TCRs and the strength of TCR-tumour antigen binding (i.e. TCR-tumour antigen binding affinity)~\citep{kastritis2013binding,stone2009t}.

The presence of tumour cells expressing non-self antigens triggers the clonal expansion of CTLs with matching TCRs. Thereupon, CTL numbers are kept under control by self-regulation mechanisms~\citep{garrod2012dissecting, yu2018recent,min2018spontaneous, stockinger2004concept, troy2003cutting}, which play a key role in the prevention of autoimmunity. Furthermore, epigenetic and genetic processes inducing stochastic and heritable changes in the antigen expression profiles of tumour cells foster dynamical intercellular variability in the expression levels of tumour antigens~\citep{campoli2008hla,sigalotti2004intratumor,urosevic2005expression}. Due to limitations posed by self-regulation mechanisms upon the numbers of CTLs targeted against different antigens at the same tumour site, such a form of intratumour heterogeneity creates the substrate for adaptation of tumour cells to the antigen-specific cytotoxic activity of CTLs and triggers adaptive changes in the repertoire of CTLs. This results in coevolutionary dynamics whereby CTLs dynamically sculpt the
antigenic distribution of tumour cells, and tumour cells concurrently reshape the repertoire of CTLs~\citep{phillips2002immunology}. 

The observation that the numbers of CD8+ and CD3+ T lymphocytes at the tumour site correlate with prognosis in different types of cancer led to the development of the `immunoscore' as a prognostic marker in cancer patients~\citep{angell2013immune, galon2006type, galon2016immunoscore, galon2019approaches}. The immunoscore provides a score that increases with the density of CD8+ and CD3+ T lymphocytes both in the centre and at the periphery of the tumour. A possible tumour classification based on the immunoscore has been proposed in~\citep{galon2006type}, where tumours with a high immunoscore are classified as `hot tumours', tumours with an intermediate immunoscore are classified as `altered tumours', and tumours with a low immunoscore are classified as `cold tumours'.
 
\paragraph{Mathematical modelling background} 
Mathematical modelling can contribute to biomedical research on the immune response to cancer by supporting experimental results with a theoretical basis, bringing new perspectives on extant empirical data, and informing new experiments~\citep{altrock2015mathematics,biselli2017organs,eftimie2011interactions,eladdadi2014mathematical,handel2020simulation,makaryan2020modeling,wilkie2013review}. In particular, different aspects of the coevolutionary dynamics between immune cells and tumour cells have been studied, under the assumption of spatially-homogeneous mixing between cells, using a number of deterministic continuum models formulated as ordinary differential equations~\citep{almuallem2021oncolytic,aguade2020tumour,balachandran2017identification,Cattani2010,d2008metamodeling,de2005validated,Frascoli2014,griffiths2020circulating, konstorum2017addressing, kuznetsov1994nonlinear,Kuznetsov2001,LinErikson2009,luksza2017neoantigen,mayer2019regulation,Takayanagi2001,walker2016concept,Wilkie2013}, integro-differential equations~\citep{bellomo2000modelling, delitala2013recognition, Kolev2003, lorenzi2015mathematical} and partial differential equations~\citep{atsou2020size, al2012evasion}.

Although more amenable to analytical and numerical approaches, which allow for a more in-depth theoretical understanding of the underlying cellular dynamics, deterministic continuum models make it more difficult to incorporate the finer details of the coevolution between tumour cells and CTLs. Moreover, they cannot easily capture adaptive phenomena that are driven by stochastic effects in the evolutionary paths of single cells. Hence, one ideally wants to derive them as the appropriate limit of stochastic discrete models (i.e. individual-based models) of the coevolutionary dynamics between tumour cells and CTLs. These individual-based models track the dynamics of single cells, thus permitting the representation of single-cell-scale mechanisms, and account for possible stochastic intercellular variability in the evolutionary trajectories of individual cells. Integrating the results of computational simulations of stochastic discrete models with analytical and numerical results of their deterministic continuum counterparts makes it possible to clearly identify the validity domain of such results, thus leading to more robust biological insight.

\paragraph{Contents of the article} In this vein, we develop an individual-based model for the coevolutionary dynamics between tumour cells and CTLs in a well-mixed system (i.e. spatial interactions are not incorporated into the model). While being aware of the fact that a variety of different cells and molecules take part in this process, here we focus on interactions involving tumour cells and CTLs only. Every cell is viewed as an individual agent whose phenotypic state is modelled by a discrete variable. For tumour cells this variable represents a parameterisation of the antigen expression profiles, while for CTLs it represents a parameterisation of TCRs.

The model takes into account the effects of the following biological processes: proliferation and death of tumour cells and CTLs; heritable, spontaneous phenotypic changes of tumour cells resulting in variation of antigenic expression; antigen-driven expansion of CTLs (i.e. {\it in situ} clonal expansion following antigen recognition); death of tumour cells due to antigen-specific cytotoxic activity of CTLs. These processes are incorporated into the model through a set of rules that correspond to a discrete-time branching random walk on the space of phenotypic states~\citep{chisholm2016evolutionary, hughes1995random}.

We show that a generalised version of the mathematical model of immune competition presented  in~\citep{lorenzi2015mathematical} can be formally obtained as the deterministic continuum limit of the individual-based model presented here. This continuum model comprises a non-local partial differential equation (PDE) for the phenotype distribution of tumour cells coupled with an integro-differential equation (IDE) for the phenotype distribution of CTLs, and shares some similarities with non-local predator-prey models such as those considered, for instance, in~\citep{delitala2013evolutionary, genieys2006adaptive,segal2013pattern,tian2017nonlocal}. In addition to the biological phenomena incorporated into the model considered in~\citep{lorenzi2015mathematical}, the deterministic continuum counterpart of the individual-based model developed here takes also into account the effect of changes in antigen expression profiles of tumour cells and more general forms of competitive feedback mechanisms regulating the growth of the numbers of tumour cells and CTLs.

The biologically relevant homogeneous steady-state solutions of the continuum model equations are found. Linear-stability analysis of these steady-state solutions is then carried out in order to identify possible conditions on the model parameters that may lead to different outcomes of immune competition and to the emergence of patterns of phenotypic coevolution between tumour cells and CTLs. We report on computational results of the individual-based model, and show that there is a good agreement between them and analytical and numerical results of the continuum model. Moreover, we explore possible scenarios in which differences between the outputs of the individual-based and continuum models may emerge. The results obtained disentangle the role of different cell parameters in the coevolutionary dynamics between tumour cells and CTLs. 

\paragraph{Structure of the article} In Section~\ref{Stochastic individual-based model for the co-evolution between CTLs and cancer cells}, the individual-based model is introduced. In Section~\ref{Formal derivation of the continuum model corresponding to the individual-based model}, the deterministic continuum counterpart of this model is presented (a formal derivation is provided in Appendix~\ref{Appendix A}). In Section~\ref{Analytical results}, the homogeneous steady-state solutions of the continuum model equations are identified and their linear stability is investigated. In Section~\ref{Computational results and critical analysis}, computational results of the individual-based model are discussed and integrated with numerical solutions of the continuum model. In Section~\ref{conclusion}, key biological implications of the main findings of this study are summarised and directions for future research are outlined. 

\section{Individual-based model}
\label{Stochastic individual-based model for the co-evolution between CTLs and cancer cells}
We model the coevolutionary dynamics between a population of tumour cells and a population of CTLs in a well-mixed system. 
The population of tumour cells is labelled by the letter $C$, while the population of CTLs is labelled by the letter $T$. Building on the modelling approach developed in~\citep{delitala2013mathematical, delitala2013recognition, lorenzi2015mathematical}, at any time $t \in [0,t_f] \subset \mathbb{R}^+$, the phenotypic state of every tumour cell is modelled by a variable $u \in \overline{\mathcal{I}}$, where $\overline{\mathcal{I}} := [-L,L] \subset \mathbb{R}$ is the closure of the set $\mathcal{I} := (-L,L) \subset \mathbb{R}$ with $L>0$, and the phenotypic state of every CTL is modelled by a variable $v \in \overline{\mathcal{I}}$. We make the following modelling assumptions.
\begin{assumption}  
\label{Ass1}
The variable $u$ represents a parameterisation of the antigen expression profiles of tumour cells, while the variable $v$ represents a parameterisation of the target antigens of the TCRs. As a result, CTLs in the phenotypic state $v$ will be primarily capable of recognising tumour cells in the phenotypic state $u=v$. 
\end{assumption}  

\begin{assumption}  
\label{Ass2}
%Tumour cells whose phenotypic states are modelled by closer values of $u$ will have a higher antigenic similarity. 
Tumour cells will have higher antigenic similarity if their phenotypic states are modelled by closer value of $u$. Hence, depending on the range of TCR affinity, CTLs in the phenotypic state $v=u$ may also be capable of recognising tumour cells in phenotypic states which are sufficiently close to $u$.
\end{assumption}  

\begin{assumption}  
\label{Ass3}
Tumour cells in similar phenotypic states (i.e. phenotypic states that are modelled by sufficiently close values of $u$) will occupy similar niches and, therefore, a form of intra-population competition (i.e. clonal competition) will occur between them. Moreover, self-regulation mechanisms act on CTLs in similar phenotypic states (i.e. phenotypic states that are modelled by sufficiently close values of $v$). Hence, a form of intra-population competition will occur between these cells as well. 
\end{assumption}  

Building upon the ideas presented in~\citep{ardavseva2020comparative, chisholm2016evolutionary, stace2019phenotype}, we model each cell as an agent that occupies a position on a lattice, which represents the space of phenotypic states. We discretise the time variable and the phenotypic states, respectively, as 
$$
t_{h} = h \tau \in [0,t_f], \quad u_{i} = i  \chi \in \overline{\mathcal{I}} \quad \text{and} \quad v_{j} = j  \chi \in \overline{\mathcal{I}}, \quad h\in\mathbb{N}_0, \;i, j \in \mathbb{Z},
$$
where $\tau \in \mathbb{R}^+_*$ and $\chi \in \mathbb{R}^+_*$ are the time- and phenotype-step, respectively. We introduce the dependent variables $N_{C_{i}}^{h}\in\mathbb{N}_0$ and $N_{T_{j}}^{h} \in\mathbb{N}_0$ to represent, respectively, the number of tumour cells on lattice site $i$ (i.e. in the phenotypic state $u_i$) and the number of CTLs on lattice site $j$ (i.e. in the phenotypic state $v_j$) at time-step $h$. The population density functions of tumour cells and CTLs (i.e. the phenotype distributions of the two cell populations) are defined, respectively, as
\begin{equation} 
n_{C_{i}}^{h} \equiv n_C(u_{i}, t_h) :=\frac{N_{C_{i}}^h}{\chi} \quad \text{and} \quad n_{T_{j}}^{h} \equiv n_T(v_{j}, t_h) :=\frac{N_{T_{j}}^h}{\chi},
\label{popdens}
\end{equation}
while the total numbers of tumour cells and CTLs (i.e. the sizes of the cell populations $C$ and $T$) are defined, respectively, as
\begin{equation} 
\rho_C^{h} \equiv \rho_C(t_h) :=\sum_i N_{C_{i}}^h \quad \text{and} \quad \rho_T^{h} \equiv \rho_T(t_h) :=\sum_j N_{T_{j}}^h.
\label{popdens_rho}
\end{equation}
In the mathematical framework of our model, the function
\begin{equation}
I_h \equiv I(t_h) := \frac{\rho_T(t_h)}{\rho_C(t_h)}
\label{immunoscore}
\end{equation}
provides a possible simplified measure of the immune score at the $h^{th}$ time-step in the well-mixed system considered here. In particular, abstracting from the immune-score based classification of tumours recalled in Section~\ref{sec:intro}, throughout the article we will refer to situations in which the average value of $I$, i.e. the quantity 
\begin{equation}
\label{def:avgI}
\overline{I} = \dfrac{\tau}{t_f} \sum_{h} I_h,
\end{equation}
is smaller than $1$ or at least about one order of magnitude larger than $1$ as `\emph{cold tumour-like scenarios}' and `\emph{hot tumour-like scenarios}', respectively, whereas the remaining situations will be classified as `\emph{altered tumour-like scenarios}'.

As mentioned earlier, we focus on a biological scenario whereby: (i) cells in the two populations divide and die due to intra-population competition (i.e. clonal competition amongst tumour cells and self-regulation of CTLs); (ii) tumour cells undergo heritable, spontaneous phenotypic changes which result in variation of antigen expression profiles; (iii) CTLs undergo antigen-driven expansion (i.e. {\it in situ} clonal expansion following antigen recognition); (iv) tumour cells die due to the antigen-specific cytotoxic activity of CTLs. These biological phenomena are incorporated into the model via the modelling strategies described in the following subsections.			

\subsection{Mathematical modelling of cell division and death due to intra-population competition}
\label{sec:modelling1}
We assume that a dividing cell is replaced by two identical cells that inherit the phenotypic state of the parent cell (i.e. the progenies are placed on the same lattice site as their parent), while a dying cell is removed from the population. \\

At every time-step $h$, we allow tumour cells and CTLs to undergo cell division at rates $\alpha_C>0$ and $\alpha_T>0$, respectively. Moreover, on the basis of Assumption~\ref{Ass3}, at every time-step $h$, we allow tumour cells in the phenotypic state $u_i$ and CTLs in the phenotypic state $v_j$ to die due to intra-population competition at rates $\mu_C \, K_{C_{i}}^{h}$ and $\mu_T \, K_{T_{j}}^{h}$, respectively, where $\mu_C, \mu_T > 0$ and
\begin{equation}
K_{C_{i}}^{h} \equiv K_{C}(u_i, t_h) := \sum_k g(u_i, u_k; \theta_C) \, N_{C_{k}}^{h}, \qquad K_{T_{j}}^{h} \equiv K_{T}(v_j, t_h) := \sum_k g(v_j, v_k; \theta_T) \, N_{T_{k}}^{h}.
\label{KCT}
\end{equation}
The function $g$ is defined as follows
\begin{equation}
\label{g}
g(x,y; \xi):=
\begin{cases}
\displaystyle{\dfrac{1}{|\mathcal{L}_{\xi}(x)|} \quad \text{if} \quad \vert y-x\vert \leq \xi} \\
\displaystyle{0  \quad\quad\quad\;\;\;\; \text{if} \quad \vert y-x\vert >\xi},
\end{cases}
\quad
\text{for} \; (x, y; \xi) \in \mathcal{I} \times \mathcal{I} \times (0, |\mathcal{I}|],
\end{equation}	
where $|\mathcal{I}|$ denotes the size of the interval $\mathcal{I}$ (i.e. $|\mathcal{I}| = 2 L$) and $|\mathcal{L}_{\xi}(x)|$ denotes the size of the interval
\begin{equation}
\label{def:Ltheta}
\mathcal{L}_{\xi}(x) := \{y \in \mathcal{I} : |y-x| \leq \xi \}, \quad (x; \xi) \in \mathcal{I} \times (0, |\mathcal{I}|].
\end{equation}	
The quantity $K_{C_{i}}^{h}$ defined via~\eqref{KCT}-\eqref{def:Ltheta} represents the number of tumour cells whose phenotypic states are modelled by values of the variable $u_k$ at a distance smaller than or equal to $\theta_C$ from $u_i$, rescaled to $|\mathcal{L}_{\theta_C}(u_i)|$. Similarly, the quantity $K_{T_{j}}^{h}$ defined via~\eqref{KCT}-\eqref{def:Ltheta} represents the number of CTLs whose phenotypic states are modelled by values of the variable $v_k$ at a distance smaller than or equal to $\theta_T$ from $v_j$, rescaled to $|\mathcal{L}_{\theta_T}(v_j)|$. Hence, the inverse of the parameter $0 < \theta_C \leq |\mathcal{I}|$ (i.e. $1/\theta_C$) provides a measure of the level of selectivity of clonal competition amongst tumour cells and the inverse of the parameter $0 < \theta_T \leq |\mathcal{I}|$ (i.e. $1/\theta_T$) provides a measure of the level of selectivity of self-regulation mechanisms acting on CTLs (i.e. the smaller $\theta_C$ and $\theta_T$, the higher the corresponding levels of selectivity). Furthermore, the parameters $\mu_C$ and $\mu_T$ represent the rates of death of tumour cells and CTLs due to these forms of intra-population competition.

\subsection{Mathematical modelling of phenotypic changes in tumour cells}
\label{sec:modelling2}
Building on the modelling strategies proposed in ~\citep{ardavseva2020comparative, chisholm2016evolutionary, stace2019phenotype}, we account for spontaneous, heritable phenotypic changes which result in variation of antigen expression profiles by allowing tumour cells to update their phenotypic states according to a random walk. More precisely, between the time-steps $h$ and $h+1$, every tumour cell enters a new phenotypic state, with probability $0 < \lambda_C < 1$, or remains in its current phenotypic state, with probability $1-\lambda_C$. When a tumour cell in the phenotypic state $u_i$ undergoes a phenotypic change, it enters into either the phenotypic state $u_{i-1}$ or the phenotypic state $u_{i+1}$ with probability $\lambda_C/2$. This models the fact that phenotypic changes occur randomly due to non-genetic instability, rather than being induced by selective pressures~\citep{huang2013genetic}. No-flux boundary conditions are implemented by aborting any attempted phenotypic variation of a tumour cell if it requires moving into a phenotypic state outside the interval $\overline{\mathcal{I}}$.	
\subsection{Mathematical modelling of tumour-immune competition}
\label{sec:modelling3}
Similarly to cell division, we assume that a CTL undergoing antigen-driven expansion is replaced by two identical cells that inherit the phenotypic state of the parent cell. Moreover, similarly to cell death due to intra-population competition, we assume that a tumour cell dying due to the antigen-specific cytotoxic activity of CTLs is removed from the population. \\

On the basis of Assumptions~\ref{Ass1} and~\ref{Ass2}, at every time-step $h$ we allow CTLs in the phenotypic state $v_j$ to undergo antigen-driven expansion at rate $\zeta_T \, \gamma \, J_{T_{j}}^{h}$, while tumour cells in the phenotypic state $u_i$ will die due to antigen-specific cytotoxic activity of CTLs at rate $\zeta_C \, \gamma \, J_{C_{i}}^{h}$. Here, $\zeta_C, \zeta_T, \gamma > 0$ and
\begin{equation}
\label{JTC}
J_{T_{j}}^{h} \equiv J_{T}(v_j, t_h) := \sum_i g(v_j,u_i; \eta) \, N_{C_{i}}^{h}, \qquad J_{C_{i}}^{h} \equiv J_{C}(u_i, t_h) := \sum_j g(u_i,v_j; \eta) \, N_{T_{j}}^{h},
\end{equation}
where the function $g$ is defined via~\eqref{g} and~\eqref{def:Ltheta}. The quantity $J_{T_{j}}^{h}$ defined via~\eqref{g}-\eqref{JTC} represents the number of tumour cells whose phenotypic states are modelled by values of the variable $u_i$ at a distance smaller than or equal to $\eta$ from $v_j$, rescaled to $|\mathcal{L}_{\eta}(v_j)|$. Similarly, the quantity $J_{C_{i}}^{h}$ defined via~\eqref{g}-\eqref{JTC} represents the number of CTLs in phenotypic states which are modelled by values of the variable $v_j$ at a distance smaller than or equal to $\eta$ from $u_i$, rescaled to $|\mathcal{L}_{\eta}(u_i)|$. Hence, the parameter $0 < \eta \leq |\mathcal{I}|$ provides a measure of the affinity range of TCRs. In more detail, this parameter determines the range of tumour antigens each CTL can
recognise: large values of $\eta$ correspond to a CTL population that is able to eliminate tumour cells expressing a large spectrum of antigens, whereas low values of $\eta$ correspond to scenarios where CTLs can only recognise tumour cells with a more specific antigenic expression. Furthermore, the parameter $\gamma$ provides a measure of the TCR-tumour antigen binding affinity (i.e. the strength of TCR-tumour antigen binding). Previous experimental studies have highlighted the role played by TCR-tumour antigen binding affinity in anti-tumour immune response and its link with immune efficiency~\citep{june2018car,gerdemann2011cytotoxic, than2015tcell}. Finally, the parameter $\zeta_T$ represents the rate at which a CTL undergoing antigen-driven expansion divides (i.e. the rate of cell division corresponding to \textit{in situ} clonal expansion), and the parameter $\zeta_C$ represents the rate at which a tumour cell can die due to the antigen-specific cytotoxic activity of a CTL. 

\subsection{Computational implementation of the individual-based model}
Numerical simulations of the individual-based model are performed in \textsc{Matlab}. At each time-step, we follow the procedures described hereafter to simulate phenotypic variation, cell division and death of tumour cells, and cell division and death of CTLs. The random numbers $r_1$, $r_2$ and $r_3$ mentioned below are real numbers drawn from the standard uniform distribution on the interval $(0,1)$ using the  built-in function \textsc{rand}.

\paragraph{Computational implementation of spontaneous, heritable phenotypic changes of tumour cells}
For each cell in population $C$, a random number, $r_1$, is generated and used to determine whether the cell undergoes a phenotypic change (i.e. $0 < r_1 < \lambda_C$) or not (i.e.  $\lambda_C \leq r_1 <  1$). If the cell undergoes a phenotypic change, then a second random number, $r_2$, is generated. If $0 < r_2 < 1/2$, then the cell moves into the phenotypic state to the left of its current state ({\it i.e.} a cell in the phenotypic state $u_i$ will move into the phenotypic state $u_{i-1}=u_{i}-\chi$), whereas if $1/2 \leq r_2 < 1$ then the cell moves into the phenotypic state to the right of its current state ({\it i.e.} a cell in the phenotypic state $u_i$ will move into the phenotypic state $u_{i+1}=u_{i}+\chi$). As mentioned earlier, no-flux boundary conditions are implemented by aborting attempted phenotypic changes that would move a cell into a phenotypic state outside the interval $\overline{\mathcal{I}}$.

\paragraph{Computational implementation of cell division and death of tumour cells and CTLs} 
For each population, the number of cells in each phenotypic state is counted. The quantities $K_C$ and $K_T$ are computed via~\eqref{KCT} and the quantities $J_C$ and $J_T$ are computed via~\eqref{JTC}, and the following definitions are used to calculate the probabilities of cell division, death and quiescence (i.e. no division nor death) for every phenotypic state of cells in populations $C$ and $T$, respectively,
    \begin{equation} \label{pC}
        {\rm P}^b_C := \tau \, \alpha_C,  \quad   {\rm P}^d_C := \tau \, \left(\mu_C \, K_{C_{i}}^{h} + \zeta_C \,\gamma\, J_{C_{i}}^{h}\right), \quad {\rm P}^q_C := 1- \left({\rm P}^b_C + {\rm P}^d_C \right)
    \end{equation}
    and
    \begin{equation} \label{pT}
        {\rm P}^b_T := \tau \, \left(\alpha_T + \zeta_T \,\gamma\, J_{T_{j}}^{h}\right),  \quad   {\rm P}^d_T := \tau \, \mu_T \, K_{T_{j}}^{h}, \quad {\rm P}^q_T := 1- \left({\rm P}^b_T + {\rm P}^d_T \right).
    \end{equation}
Notice that we are implicitly assuming that the time-step $\tau$ is sufficiently small that $0 < {\rm P}^k_C < 1$ and $0 < {\rm P}^k_T < 1$ for all $k \in \{b, d, q\}$. For each cell, a random number, $r_3$, is generated and the cells' fate is determined by comparing this number with the probabilities of division, death and quiescence corresponding to the phenotypic state of the cell. In more detail, for a cell in population $C$: if $0 < r_3 < {\rm P}^d_C$ then the cell is considered dead and is removed from the population; if  ${\rm P}^d_C \leq r_3 < {\rm P}^d_C + {\rm P}^b_C$ then the cell undergoes division and an identical daughter cell is created; whereas if ${\rm P}^d_C + {\rm P}^b_C \leq r_3 < 1$ then the cell remains quiescent (i.e. does not divide nor die). Similarly, for a cell in population $T$: if $0 < r_3 < {\rm P}^d_T$ then the cell is considered dead and is removed from the population; if  ${\rm P}^d_T \leq r_3 < {\rm P}^d_T + {\rm P}^b_T$ then the cell undergoes division and an identical daughter cell is created; whereas if ${\rm P}^d_T + {\rm P}^b_T \leq r_3 < 1$ then the cell remains quiescent. 

\section{Corresponding deterministic continuum model}	
\label{Formal derivation of the continuum model corresponding to the individual-based model}
In the case where cell dynamics are governed by the rules described in Sections~\ref{sec:modelling1}-\ref{sec:modelling3}, between time-steps $h$ and $h+1$ a tumour cell in the phenotypic state $u_i$ may divide, die or remain quiescent (i.e. not divide nor die) with probabilities defined via~\eqref{pC}, while a CTL in the phenotypic state $v_j$ may divide, die or remain quiescent with probabilities defined via~\eqref{pT}. Hence, recalling that between time-steps $h$ and $h+1$ a tumour cell in the phenotypic state $u_i$ may also enter into either of the phenotypic states $u_{i-1}$ and $u_{i+1}$ with probabilities $\lambda_C/2$, the principle of mass balance gives the following system of coupled difference equations for the population densities $n_{C_{i}}^{h}$ and $n_{T_{j}}^{h}$:
\begin{equation}
\label{eq:mastereqs}
\begin{cases}
n_{C_{i}}^{h+1}=\left(2 \, {\rm P}^b_C + {\rm P}^q_C \right)\left[\dfrac{\lambda_C}{2} \left(n_{C_{i+1}}^{h} + n_{C_{i-1}}^{h}\right) + (1-\lambda_C) \, n_{C_{i}}^{h}\right], \quad (u_i, t_h) \in \mathcal{I} \times (0,t_f],
\\\\
n_{T_{j}}^{h+1}=\left(2 \, {\rm P}^b_T + {\rm P}^q_T \right) \, n_{T_{i}}^{h}, \quad (v_j, t_h) \in \overline{\mathcal{I}} \times (0,t_f].
\end{cases}
\end{equation} 
The difference equation~\eqref{eq:mastereqs}$_1$ for $n^h_{C_{i}}$ is subject to no-flux boundary conditions due to the fact that, as mentioned in Section~\ref{sec:modelling2}, any attempted phenotypic variation of a tumour cell is aborted if it requires moving into a phenotypic state outside the interval $\overline{\mathcal{I}}$.

Starting from the system of coupled difference equations~\eqref{eq:mastereqs}, letting the time-step $\tau \to 0$ and the phenotype-step $\chi \to 0$ in such a way that
\begin{equation}
\lambda_C\frac{\chi^2}{2\tau} \to \beta_C \quad \text{with} \quad 0 < \beta_C < \infty,
\label{conditionBeta}
\end{equation}
where the parameter $\beta_C$ is the rate of spontaneous, heritable phenotypic changes of tumour cells, using the method employed in~\citep{ardavseva2020comparative, chisholm2016evolutionary, stace2019phenotype}, it is possible to formally show (see Appendix~\ref{Appendix A}) that the deterministic continuum counterpart of the stochastic discrete model comprises the following PDE-IDE system for the cell population density functions $n_{C}(u,t)$ and $n_{T}(v,t)$
\begin{equation}
\label{eq:PDEs}
	    	 \begin{dcases*} \partial_t n_C - \beta_C \, \partial^2_{uu} n_C = \Big[\alpha_C - \mu_C \, K_C(u,t) - \zeta_C \, \gamma \, J_C(u,t)\Big] n_C, \quad (u,t) \in \mathcal{I}\times (0,t_f], \\
		 \partial_t n_T = \Big[\alpha_T - \mu_T K_T(v,t) + \zeta_T \, \gamma \, J_T(v,t) \Big] n_T, \quad (v,t) \in \overline{\mathcal{I}}\times (0,t_f],
		 \\
	    			J_C(u,t):=\int_{\mathcal{I}} g(u,v; \eta) \, n_T(v,t) \, \mathrm{d}v, \quad K_C(u,t):=\int_{\mathcal{I}} g(u,w; \theta_C) \, n_C(w,t) \, \mathrm{d}w,  
				\\
		J_T(v,t):=\int_{\mathcal{I}} g(v,u; \eta) \, n_C(u,t)\mathrm{d}u, \quad K_T(v,t):=\int_{\mathcal{I}} g(v,w; \theta_T) \, n_T(w,t) \, \mathrm{d}w,
\end{dcases*}
\end{equation}
with $\mathcal{I}=(-L,L)$. Here, the function $g$ is defined via~\eqref{g} and~\eqref{def:Ltheta}, and the non-local PDE~\eqref{eq:PDEs}$_1$ for $n_C$ is subject to the following no-flux boundary conditions
\begin{equation}
   \partial_u n_C(-L,t)=0 \quad \text{and} \quad \partial_u n_C(L,t)=0 \quad \text{for all } t \in (0,t_f].
    \label{neumann}
\end{equation}
We remark that linear diffusion operators like the one in the PDE~\eqref{eq:PDEs}$_1$, which are the deterministic, continuum counterparts of underlying random walks over the space of phenotypic states, have been widely used to model the effect of heritable, spontaneous phenotypic changes in cell populations -- see, for instance, the review article~\citep{ChisholmBBAGS2016} and references therein.
\section{Steady-state and linear-stability analyses of the  continuum model equations}
\label{Analytical results}
In this section, we first identify the biologically relevant homogeneous steady-state solutions of the continuum model equations. Then, we carry out linear-stability analysis to: (i) determine conditions that may lead to the eradication of tumour cells by CTLs or to the coexistence between the two cell populations, and (ii) identify sufficient conditions for the emergence of patterns of phenotypic coevolution between tumour cells and CTLs.

\subsection{Biologically relevant steady-state solutions} 
A biologically relevant steady-state solution of the PDE-IDE system~\eqref{eq:PDEs} subject to the boundary conditions~~\eqref{neumann} is given by a pair of real, non-negative functions $n_C^*(u)$ and $n_T^*(v)$ that satisfy the following system
\begin{equation}
\label{eqsteady2}
	    	 \begin{dcases*} - \beta_C \, \partial^2_{uu} n_C^* = \Big[\alpha_C - \mu_C \, K_C^*(u) - \zeta_C \, \gamma \, J_C^*(u)\Big] n_C^*, \quad u \in \mathcal{I}, \\
		 \Big[\alpha_T - \mu_T K_T^*(v) + \zeta_T \, \gamma \, J_T^*(v) \Big] n_T^*=0, \quad v \in \overline{\mathcal{I}},
\end{dcases*}
\end{equation}
where $\mathcal{I}=(-L,L)$, with \eqref{eqsteady2}$_1$ subject to the boundary conditions
\begin{equation}
   \partial_u n_C^*(-L)=0 \quad \text{and} \quad \partial_u n_C^*(L)=0.
    \label{neumannss}
\end{equation}
In the system~\eqref{eqsteady2},
\begin{equation}
\label{eq:JCs}
J_C^*(u):=\int_{\mathcal{I}} g(u,v; \eta) \, n_T^*(v) \, \mathrm{d}v, \quad K_C^*(u):=\int_{\mathcal{I}} g(u,w; \theta_C) \, n_C^*(w) \, \mathrm{d}w
\end{equation}
and
\begin{equation}
\label{eq:JTs}
		J_T^*(v):=\int_{\mathcal{I}} g(v,u; \eta) \, n_C^*(u) \, \mathrm{d}u, \quad K_T^*(v):=\int_{\mathcal{I}} g(v,w; \theta_T) \, n_T^*(w) \, \mathrm{d}w.
\end{equation}
The components of homogeneous steady-state solutions satisfy the following system of equations
\begin{equation}
\label{eqsteady2hss}
	    	 \begin{dcases*} \Big[\alpha_C - \mu_C \, K_C^*(u) - \zeta_C \, \gamma \, J_C^*(u)\Big] n_C^* = 0 , \quad u \in \overline{\mathcal{I}}, \\
		 \Big[\alpha_T - \mu_T K_T^*(v) + \zeta_T \, \gamma \, J_T^*(v) \Big] n_T^*=0, \quad v \in \overline{\mathcal{I}}
\end{dcases*}
\end{equation}
and are of the form
\begin{equation}\label{eq:homsteady}
n_C^*(u) = \dfrac{\rho_C^*}{|\mathcal{I}|}  \;\; \forall u \in \overline{\mathcal{I}} \quad \text{and} \quad n_T^*(v) = \dfrac{\rho_T^*}{|\mathcal{I}|}  \;\; \forall v \in \overline{\mathcal{I}},
\end{equation}
where $\rho_C^* \geq 0$ and $\rho_T^* \geq 0$ satisfy the following system of algebraic equations
  	   						\begin{equation}
  	   			\begin{dcases*}
	    	 \Big(\alpha_C \, |\mathcal{I}| -\mu_C \, \rho_C^*-\gamma_C \, \rho_T^*\Big)\rho_C^*=0, \\
	    	 \Big(\alpha_T \, |\mathcal{I}| -\mu_T \, \rho_T^*+ \gamma_T \, \rho_C^*\Big)\rho_T^*=0,
	    	  \end{dcases*} 
		  \quad
		  \text{with}
		  \quad
		  \gamma_C:=\zeta_C \, \gamma \quad \text{and} \quad \gamma_T:= \zeta_T \,\gamma.
	    	  \label{ptEqrho}
  	   			\end{equation}	
The system of algebraic equations~\eqref{ptEqrho} is obtained by first integrating the PDE~\eqref{eq:PDEs}$_1$ over $\mathcal{I}$ and imposing the boundary conditions~\eqref{neumann}, then integrating the IDE~\eqref{eq:PDEs}$_2$ over $\mathcal{I}$, next substituting ansatz~\eqref{eq:homsteady} into the resulting equations and equating to zero their right-hand sides, and finally using the fact that, when the function $g$ is defined via~\eqref{g} and~\eqref{def:Ltheta}, 
\begin{equation}
\label{eq:idg}
\int_{\mathcal{I}} g(x,y; \xi) \, \mathrm{d}y = 1, \;\; \forall x \in \mathcal{I}, \;\xi\in  (0, |\mathcal{I}|].
\end{equation}
In particular, since we are studying tumour-immune competition, we are interested in solutions of the system of equations~\eqref{ptEqrho} whose $\rho_T^*$ component is strictly positive. There exist two solutions of this type, that is, the semitrivial solution
\begin{equation}
(\rho_{C1}^*,\rho_{T1}^*) =\left(0,\dfrac{\lvert\mathcal{I}\rvert\alpha_T}{\mu_T}\right), 
\label{steady_b}
\end{equation}  
and, provided that the following condition on the model parameters is met
\begin{equation}
\gamma < \dfrac{\mu_T}{\alpha_T} \dfrac{\alpha_C}{\zeta_C},
 \label{condition}
\end{equation}
the nontrivial solution
\begin{equation}(\rho_{C2}^*,\rho_{T2}^*)=\left(\lvert\mathcal{I}\rvert \, \dfrac{(\alpha_C \mu_T-\alpha_T\gamma_C)}{\gamma_T\gamma_C + \mu_C\mu_T}, \lvert\mathcal{I}\rvert \, \dfrac{(\alpha_T\mu_C+\alpha_C\gamma_T)}{\gamma_T\gamma_C + \mu_C\mu_T}\right).
\label{steady_d}
\end{equation}
The semitrivial steady-state solution given by~\eqref{eq:homsteady} and \eqref{steady_b} corresponds to biological scenarios whereby tumour cells are eradicated by CTLs, while the nontrivial steady-state solution given by~\eqref{eq:homsteady} and \eqref{steady_d} corresponds to situations where coexistence between tumour cells and CTLs occurs. Notice that condition~\eqref{condition} indicates that lower TCR-tumour antigen binding affinity (i.e. smaller values of $\gamma$) make it more likely that tumour cells survive the cytotoxic action of CTLs, thus promoting coexistence between the two cell populations.

\subsection{Linear-stability analysis} 
Linearising the PDE-IDE system~\eqref{eq:PDEs}, subject to the boundary conditions~\eqref{neumann}, about a steady-state of components $n_C^*(u)$ and $n_T^*(v)$, and using the conditions given by equations~\eqref{eqsteady2}, we obtain the following PDE-IDE system for the perturbations $\tilde{n}_{C}(u,t)$ and $\tilde{n}_{T}(v,t)$
\begin{equation}
\small
\label{eq5}
	    	 \begin{dcases*} \partial_t \tilde{n}_C - \beta_C \, \partial^2_{uu} \tilde{n}_C = \left[\alpha_C - \mu_C K_C^*(u) - \gamma_C J_C^*(u) \right] \tilde{n}_C - \Big[\mu_C \, \tilde{K}_C(u,t) + \gamma_C \, \tilde{J}_C(u,t)\Big] n_C^*, \quad (u,t) \in \mathcal{I}\times (0,t_f], \\
		 \partial_t \tilde{n}_T = \left[\alpha_T - \mu_T K_T^*(v)  + \gamma_T J_T^*(v) \right] \tilde{n}_T - \Big[\mu_T \, \tilde{K}_T(v,t) - \gamma_T \, \tilde{J}_T(v,t)\Big] n_T^*, \quad (v,t) \in \overline{\mathcal{I}}\times (0,t_f],
\end{dcases*}
\end{equation}
subject to the boundary conditions
\begin{equation}
   \partial_u \tilde n_C(-L,t)=0 \quad \text{and} \quad \partial_u \tilde n_C(L,t)=0 \quad \text{for all } t \in (0,t_f].
    \label{neumannpert}
\end{equation}
In the system~\eqref{eq5}, $J_C^*(u)$ and $K_C^*(u)$ are defined via~\eqref{eq:JCs}, $J_T^*(v)$ and $K_T^*(u)$ are defined via~\eqref{eq:JTs}, and
\begin{equation}
\label{eq:JCt}
\tilde J_C(u,t):=\int_{\mathcal{I}} g(u,v; \eta) \, \tilde n_T(v,t) \, \mathrm{d}v, \quad \tilde K_C(u,t):=\int_{\mathcal{I}} g(u,w; \theta_C) \, \tilde n_C(w,t) \, \mathrm{d}w,
\end{equation}
\begin{equation}
\label{eq:JTt}
\tilde J_T(v,t):=\int_{\mathcal{I}} g(v,u; \eta) \, \tilde n_C(u,t) \, \mathrm{d}u, \quad \tilde K_T(v,t):=\int_{\mathcal{I}} g(v,w; \theta_T) \, \tilde n_T(w,t) \, \mathrm{d}w.
\end{equation}
Due to~\eqref{eqsteady2hss}, if the steady-state solution $(n_C^*,n_T^*)$ is given by~\eqref{eq:homsteady} and~\eqref{steady_b} then the PDE-IDE system~\eqref{eq5} reduces to
\begin{equation}
\label{eq5sst1}
\begin{dcases*} \partial_t \tilde{n}_C - \beta_C \, \partial^2_{uu} \tilde{n}_C =  \Big[\alpha_C - \gamma_C J_C^*(u) \Big] \tilde{n}_C, \quad (u,t) \in \mathcal{I}\times (0,t_f], \\
\partial_t \tilde{n}_T = - \Big[\mu_T \, \tilde{K}_T(v,t) - \gamma_T \, \tilde{J}_T(v,t)\Big] \dfrac{\rho_{T1}^*}{|\mathcal{I}|}, \quad (v,t) \in \overline{\mathcal{I}}\times (0,t_f],
\end{dcases*}
\end{equation}
whereas if condition~\eqref{condition} is met and the steady-state solution $(n_C^*,n_T^*)$ is given by~\eqref{eq:homsteady} and~\eqref{steady_d} then the PDE-IDE system~\eqref{eq5} reduces to
\begin{equation}
\label{eq5sst2}
	    	 \begin{dcases*} \partial_t \tilde{n}_C - \beta_C \, \partial^2_{uu} \tilde{n}_C = - \Big[\mu_C \, \tilde{K}_C(u,t) + \gamma_C \, \tilde{J}_C(u,t)\Big]\dfrac{\rho_{C2}^*}{|\mathcal{I}|}, \quad (u,t) \in \mathcal{I}\times (0,t_f], \\
		 \partial_t \tilde{n}_T = - \Big[\mu_T \, \tilde{K}_T(v,t) - \gamma_T \, \tilde{J}_T(v,t)\Big] \dfrac{\rho_{T2}^*}{|\mathcal{I}|}, \quad (v,t) \in \overline{\mathcal{I}}\times (0,t_f].
\end{dcases*}
\end{equation}

\subsubsection{Conditions for eradication of tumour cells by CTLs or coexistence between the two cell populations} 
\label{temporal perturbations}
In order to determine conditions on the model parameters that may lead to the eradication of tumour cells by CTLs or to the coexistence between the two cell populations, we study the stability of the steady-state solutions given by~\eqref{eq:homsteady} and~\eqref{steady_b} or~\eqref{steady_d} to perturbations of the form
\begin{equation}
\tilde{n}_C(u,t) = \epsilon_C \, e^{\lambda t}  \;\; \forall u \in \overline{\mathcal{I}} \quad \text{ and } \quad \tilde n_T(v,t) = \epsilon_T \, e^{\lambda t}  \;\; \forall v \in \overline{\mathcal{I}} \quad \text{with} \quad \epsilon_C, \epsilon_T \in \mathbb{R}_*, \; \lambda \in \mathbb{C}.
\label{timedeppert}
\end{equation}

Substituting the ansatz~\eqref{timedeppert} into the PDE-IDE system~\eqref{eq5sst1} and using property~\eqref{eq:idg} along with the expression~\eqref{steady_b} of $\rho_{T1}^*$ gives the following system of algebraic equations
\begin{equation}
\label{eq10a}
	    	 \begin{dcases*} \lambda \epsilon_C = \Big(\alpha_C-\gamma_C \, \dfrac{\alpha_T}{\mu_T}\Big)\epsilon_C,  \\
	    			\lambda \epsilon_T = - \Big(\mu_T \,\epsilon_T - \gamma_T \, \epsilon_C\Big) \, \dfrac{\alpha_T}{\mu_T},
	    			 \end{dcases*}
\end{equation}
which can be written in matrix form as 
	    	\[
\begin{bmatrix}
   \alpha_C -\gamma_C \dfrac{\alpha_T}{\mu_T} -\lambda       &  0 \\
   \gamma_T \dfrac{\alpha_T}{\mu_T} & - \alpha_T -\lambda
\end{bmatrix}
\begin{bmatrix}
\epsilon_C \\
\epsilon_T
\end{bmatrix}
=
\begin{bmatrix}
   0\\
   0
\end{bmatrix}.
\]
For a non-trivial solution of system~\eqref{eq10a} to exist, the determinant of the above matrix must be zero. This leads to the following quadratic equation for $\lambda$
$$ 
\lambda^2-B\lambda+C=0
$$
with
$$
B := \alpha_C - \gamma_C \dfrac{\alpha_T}{\mu_T} - \alpha_T \quad \text{and} \quad C := \alpha_T\left(\gamma_C \frac{\alpha_T}{\mu_T} -\alpha_C\right).
$$
Hence, the semitrivial steady-state solution given by~\eqref{eq:homsteady} and~\eqref{steady_b} is locally asymptotically stable if the reverse of condition~\eqref{condition} holds, that is if
\begin{equation}
\gamma > \dfrac{\mu_T}{\alpha_T} \dfrac{\alpha_C}{\zeta_C},
\label{condition1}
\end{equation}
since in this case $B<0$ and $C>0$ (i.e. ${\rm Re}(\lambda)<0$).
On the other hand, performing similar calculations on the PDE-IDE system~\eqref{eq5sst2} gives the the following system of algebraic equations
\begin{equation}
\label{eq10an1}
	    	 \begin{dcases*} \lambda \epsilon_C = - \Big(\mu_C \, \epsilon_C+\gamma_C \, \epsilon_T\Big) \, \dfrac{(\alpha_C \mu_T-\alpha_T\gamma_C)}{\gamma_T\gamma_C + \mu_C\mu_T},  \\
	    			\lambda \epsilon_T = - \Big(\mu_T \,\epsilon_T - \gamma_T \, \epsilon_C\Big) \, \dfrac{(\alpha_T\mu_C+\alpha_C\gamma_T)}{\gamma_T\gamma_C + \mu_C\mu_T},
	    			 \end{dcases*}
\end{equation}
which can be written in matrix form as 
	    	\[
\begin{bmatrix}
    - \mu_C \dfrac{(\alpha_C \mu_T-\alpha_T\gamma_C)}{\gamma_T\gamma_C + \mu_C\mu_T} -\lambda       &  -\gamma_C \dfrac{(\alpha_C \mu_T-\alpha_T\gamma_C)}{\gamma_T\gamma_C + \mu_C\mu_T} \\
   \gamma_T \dfrac{(\alpha_T\mu_C+\alpha_C\gamma_T)}{\gamma_T\gamma_C + \mu_C\mu_T} &  - \mu_T \dfrac{(\alpha_T\mu_C+\alpha_C\gamma_T)}{\gamma_T\gamma_C + \mu_C\mu_T} -\lambda
\end{bmatrix}
\begin{bmatrix}
\epsilon_C \\
\epsilon_T
\end{bmatrix}
=
\begin{bmatrix}
   0\\
   0
\end{bmatrix}.
\]
For a non-trivial solution of system~\eqref{eq10an1} to exist, the determinant of the above matrix must be zero. This leads to the following quadratic equation for $\lambda$
$$ 
\lambda^2-B\lambda+C=0
$$
with
$$
B := - \left[\mu_C \dfrac{(\alpha_C \mu_T-\alpha_T\gamma_C)}{\gamma_T\gamma_C + \mu_C\mu_T} + \mu_T \dfrac{(\alpha_T\mu_C+\alpha_C\gamma_T)}{\gamma_T\gamma_C + \mu_C\mu_T} \right]\quad $$ \text{and} $$ \quad C := \dfrac{(\alpha_C \mu_T-\alpha_T\gamma_C) (\alpha_T\mu_C+\alpha_C\gamma_T)}{(\gamma_T\gamma_C + \mu_C\mu_T)^2}\left[\mu_C\mu_T+\gamma_C\gamma_T\right].
$$
Hence, if condition~\eqref{condition} is met, then the nontrivial steady-state solution given by~\eqref{eq:homsteady} and~\eqref{steady_d} is locally asymptotically stable, since in this case $B<0$ and $C>0$ (i.e. ${\rm Re}(\lambda)<0$).

\subsubsection{Conditions for the emergence of patterns of phenotypic coevolution between tumour cells and CTLs} 
\label{spatio-temporal perturbations}
In order to identify sufficient conditions for the emergence of patterns of phenotypic coevolution between tumour cells and CTLs, we study the stability of the nontrivial steady-state solution given by~\eqref{eq:homsteady} and~\eqref{steady_d} to perturbations of the form
\begin{equation}
\tilde{n}_C(u,t) = \epsilon_C \, e^{\lambda t} \, \varphi_k(u) \quad \text{ and } \quad \tilde n_T(v,t) = \epsilon_T \, e^{\lambda t} \, \varphi_k(v) \quad \text{with} \quad \epsilon_C, \epsilon_T \in \mathbb{R}_*, \; \lambda \in \mathbb{C}.
\label{ansatz3}
\end{equation}
Here, $\{\varphi_k\}_{k \geq 1}$ are the eigenfunctions of the Laplace operator on $\mathcal{I}$ with homogeneous Neumann boundary conditions indexed by the wavenumber $k$, that is,
\begin{equation}
\varphi_k(x)=\cos\left(k \, x\right)\quad \text{with} \quad k=\frac{m\pi}{\lvert\mathcal{I}\rvert}, \quad m\in \mathbb{N}, \quad x \in \overline{\mathcal{I}}.
 \label{cos}
\end{equation}

Substituting the ansatz given by~\eqref{ansatz3} and~\eqref{cos} into the PDE-IDE system~\eqref{eq5sst2}, using the fact that
$$
\int_{\mathcal{I}} g(x,y; \xi) \varphi_k(y) \; {\rm d}y = \dfrac{\sin(k \xi)}{k} \, \Psi(x;\xi) \, \varphi_k(x) \quad \text{with} \quad \Psi(x;\xi):=\dfrac{2}{|\mathcal{L}_{\xi}(x)|},
$$
for all $x \in \mathcal{I}$ and $\xi\in  (0, |\mathcal{I}|]$,
we obtain the following infinite system of algebraic equations 
\begin{equation}
	    	 \begin{dcases*} \lambda \, \epsilon_C  = - k^2 \, \beta_C \, \epsilon_C - \Big(\mu_C \, \frac{\sin(k\theta_C)}{k } \, \Psi(x;\theta_C)\, \epsilon_C + \gamma_C \, \frac{\sin(k\eta)}{k }\, \Psi(x;\eta) \, \epsilon_T \Big) \, \dfrac{\rho_{C2}^*}{|\mathcal{I}|} ,\\
		 \lambda \, \epsilon_T = -\Big(\mu_T \, \frac{\sin(k\theta_T)}{k } \, \Psi(x;\theta_T)\, \epsilon_T - \gamma_T \, \frac{\sin(k\eta)}{k} \, \Psi(x;\eta)\, \epsilon_C \Big) \, \dfrac{\rho_{T2}^*}{|\mathcal{I}|},
		\end{dcases*}
	    			 \label{eq7}
\end{equation}
which can be written in matrix form as
	    	\[
\begin{bmatrix}
    - k^2 \, \beta_C - \mu_C\dfrac{\sin(k\theta_C)}{k } \, \Psi(x;\theta_C) \dfrac{\rho_{C2}^*}{|\mathcal{I}|}-\lambda       &  - \gamma_C\dfrac{\sin(k\eta)}{k }\, \Psi(x;\eta) \dfrac{\rho_{C2}^*}{|\mathcal{I}|} \\
    \gamma_T\dfrac{\sin(k\eta)}{k }\, \Psi(x;\eta) \dfrac{\rho_{T2}^*}{|\mathcal{I}|} & - \mu_T\dfrac{\sin(k\theta_T)}{k } \, \Psi(x;\theta_T) \dfrac{\rho_{T2}^*}{|\mathcal{I}|} -\lambda
\end{bmatrix}
\begin{bmatrix}
\epsilon_C \\
\epsilon_T 
\end{bmatrix}
=
\begin{bmatrix}
   0\\
   0
\end{bmatrix}.
\]
For each $x \in \mathcal{I}$, for a non-trivial solution of the system of algebraic equations~\eqref{eq7} to exist the determinant of the above matrix must be zero. For each $x \in \mathcal{I}$, this leads to the following quadratic equation for $\lambda$
$$ 
\lambda^2-B\lambda+C=0
$$
where
\begin{equation}
 B\equiv B(k,x):=-k^2 \, \beta_C -\mu_C\dfrac{\sin(k\theta_C)}{k } \dfrac{\rho_{C2}^*}{|\mathcal{I}|}\, \Psi(x;\theta_C) -  \mu_T\dfrac{\sin(k\theta_T)}{k } \dfrac{\rho_{T2}^*}{|\mathcal{I}|}\, \Psi(x;\theta_T)
\label{eq:B(k)}
\end{equation}
and 
\begin{equation*}
\begin{aligned}
C\equiv C(k,x):&= k^2 \, \beta_C \, \mu_T\dfrac{\sin(k\theta_T)}{k } \dfrac{\rho_{T2}^*}{|\mathcal{I}|}\Psi(x;\theta_T) \\
&+ \dfrac{\rho_{C2}^* \, \rho_{T2}^*}{|\mathcal{I}|^2} \left[\gamma_C \gamma_T \left(\dfrac{\sin(k\eta)}{k }\Psi(x;\eta)\right)^2 + \mu_C \mu_T \dfrac{\sin(k\theta_C)}{k } \, \dfrac{\sin(k\theta_T)}{k} \Psi(x;\theta_C)\Psi(x;\theta_T) \right].
\end{aligned}
\end{equation*}
A sufficient condition for the nontrivial steady-state solution given by~\eqref{eq:homsteady} and~\eqref{steady_d} to be driven unstable by perturbations of the form~\eqref{ansatz3} (i.e. for patterns of phenotypic coevolution between tumour cells and CTLs to be formed) is that $B>0$ and/or $C<0$ so that ${\rm Re}(\lambda)>0$ for all $x \in \mathcal{I}$. In particular, in the case where 
\begin{equation}
\label{cond:thetaeta2}
\theta_C=\theta_T=\theta, 
\end{equation}
since $k$ is defined via~\eqref{cos}, for the condition $B(k, x)>0$ to hold for all $x \in \mathcal{I}$ it suffices that
$$
\beta_C < \dfrac{1}{|\mathcal{I}|}\min_{k \in \mathcal{K}}\left\lbrace - \frac{\sin(k\theta)}{k } \left(\frac{\rho_{C2}^*\mu_C+\rho_{T2}^*\mu_T}{k^2}\right) \right\rbrace\min_{x \in \mathcal{I}} \Psi(x;\theta),
$$
where $\displaystyle{\mathcal{K} := \left\{k=\frac{m\pi}{\lvert\mathcal{I}\rvert}, \, m\in \mathbb{N} : \sin(k \theta) < 0 \right\}}$. Since, under definition~\eqref{def:Ltheta}, 
$$
\displaystyle{\min_{x \in \mathcal{I}} \Psi(x;\theta)= \dfrac{2}{\underset{x \in \mathcal{I}}{\max}|\mathcal{L}_{\theta}(x)|}=\dfrac{1}{\theta}},
$$
the above condition on $\beta_C$ reduces to 
\begin{equation}
\beta_C < \dfrac{1}{|\mathcal{I}|}\min_{k \in \mathcal{K}}\left\lbrace - \frac{\sin(k\theta)}{k \theta} \left(\frac{\rho_{C2}^*\mu_C+\rho_{T2}^*\mu_T}{k^2}\right) \right\rbrace.
\label{cond:patterns}
\end{equation}
%\begin{remark}
%\label{patternsattheboundary}
%Notice that in the case where $\theta_C=\theta_T=2 L$ and $\eta \in (0,L)$, since
%$$
%\int_{-L}^L g(x,y; 2L)  \varphi_k(x) \; {\rm d}x = \dfrac{1}{2L} \int_{-L}^L \varphi_k(x) \; {\rm d}x  = 0 \quad \forall x \in [-L,L],
%$$
%similar calculations give
%$$
%B\equiv B(k):=-k^2 \, \beta_C, \quad C\equiv C(k):= \dfrac{\rho_{C2}^* \, \rho_{T2}^*}{|\mathcal{I}|^2} \gamma_C \gamma_T \left(\dfrac{\sin(k\eta)}{k \eta}\right)^2.
%$$
%Hence, in this case, if patterns of phenotypic coevolution between tumour cells and CTLs are formed they will be initiated by the growth of perturbations in the regions $[-L, -L + \eta)$ and $(L-\eta, L]$.
%\end{remark}

\section{Numerical simulations}
\label{Computational results and critical analysis}
In this section, we report on computational results of the individual-based model along with numerical solutions of the corresponding continuum model given by the PDE-IDE system~\eqref{eq:PDEs} and subject to the boundary conditions~\eqref{neumann}. Simulations are integrated with the results of steady-state and linear-stability analyses of the continuum model equations presented in Section~\ref{Analytical results}. In particular, we investigate the way in which the outputs of the models are affected by key parameters whose impact on the coevolutionary dynamics between tumour cells and CTLs is of particular biological interest. Such key parameters are: the TCR-tumour antigen binding affinity, $\gamma$, the level of selectivity of self-regulation mechanisms acting on CTLs, $1/\theta_T$, the level of selectivity of clonal competition amongst tumour cells, $1/\theta_C$, and the affinity range of TCRs, $\eta$. Moreover, we explore the existence of scenarios in which differences between the outputs produced by the two models can emerge due to effects which reduce the quality of the approximation of the individual-based model provided by the continuum model.
\begin{comment}
In particular, we show how the model outcomes are affected by the TCR-tumour antigen binding affinity, $\gamma$, the affinity range of TCRs, $\eta$, the level of self-regulation of CTLs, $\theta_T$, and the selectivity of clonal competition amongst tumour cells, $\theta_C$.
\end{comment}
\subsection{Set-up of numerical simulations}
Without loss of generality we choose $L=1$, so that $\overline{\mathcal{I}}=[-1,1]$ and $|\mathcal{I}|=2$, and consider a discretisation of the interval $[-1,1]$ consisting of $1500$ points (i.e. the phenotype-step is $\chi\approx0.0013$). Furthermore, we use the time-step $\tau=0.05$ and, unless otherwise specified, we choose the final time $t_f=30$ days.

Building on the results of steady-state and linear-stability analyses of the continuum model equations presented in Section~\ref{Analytical results}, we carry out simulations using the following initial condition for the individual-based model
\begin{equation}
n_{C}^0(u_i) :=10^{4}(1 + a \, \cos(A\, u_i)), \qquad n_{T}^0(v_j) :=10^{4}(2 + a \, \cos(A \,v_j)), \quad a \geq 0, \; A > 0.
\label{initcondAgentbased}
\end{equation}
In Appendix~\ref{Appendix3}, we provide a detailed description of the methods employed to numerically solve the PDE-IDE system~\eqref{eq:PDEs} complemented with the boundary conditions~\eqref{neumann} and the continuum analogue of the initial condition~\eqref{initcondAgentbased}, i.e. the initial condition
\begin{equation}
n_{C}^0(u) :=10^{4}(1 + a \, \cos(A\, u)), \qquad n_{T}^0(v) :=10^{4}(2 + a \, \cos(A \,v)), \quad a \geq 0, \; A > 0.
\label{initcondPDE}
\end{equation}

Unless otherwise specified, we use the parameter values listed in Table~\ref{Tab:Tcr}. Here, the values of the parameters $\alpha_C$, $\alpha_T$, $\zeta_C$ and $\zeta_T$ are consistent with previous measurement and estimation studies on the dynamics of tumour cells and CTLs~\citep{christophe2015biased, de2009mathematical, kuznetsov1994nonlinear, schlesinger2014coevolutionary}. The values of the parameters $\mu_C$ and $\mu_T$ and the range of values of the parameters $\theta_C$ and $\theta_T$ are chosen so as to ensure that the equilibrium sizes and phenotype distributions of the two cell populations in isolation are biologically relevant. The range of values of the parameter $\eta$ is consistent with experimental estimations of the precursor frequency of CTLs~\citep{blattman2002estimating}, while the values of the parameter $\gamma$ are consistent with those used in~\citep{stromberg2006robustness,stromberg2013diversity}. The value of the parameter $\lambda_C$ is taken from~\citep{stace2019phenotype} and corresponds to values of $\beta_C$ that are consistent with experimental data reported in~\citep{doerfler2006dna,duesberg2000explaining}.

\begin{table}[ht]
	\footnotesize
\caption{Parameter values used in numerical simulations and their sources}
\centering
 	\begin{tabular}{llll}
  \hline
   & Biological meaning & Value & Source \\
  \hline
   $\alpha_C$ & Rate of tumour cell proliferation  & 1.5/day &  \citep{christophe2015biased}\\
  $\alpha_T$ & Rate of antigen-independent CTL proliferation  & 5 $\times 10^{-2}$/day & \citep{de2009mathematical}\\
  $\mu_C$ & Rate of death of tumour cells due to clonal competition& 1.5 $\times 10^{-6} \mu l$/day & \textit{ad hoc} \\
  $\mu_T$ & Rate of death of CTLs due to self-regulation mechanisms& 5 $\times 10^{-6}\mu l$/day & \textit{ad hoc}\\
  $\zeta_C$ & Killing rate of tumour cells by CTLs & 5 $\times 10^{-6}\mu l$/day  & \citep{kuznetsov1994nonlinear}\\
  $\zeta_T$ & Rate of replication of CTLs following recognition & 3 $\times 10^{-5}\mu l$/day & \citep{schlesinger2014coevolutionary}\\
  $\eta$ & Affinity range of TCRs & [0.1, 2] & \citep{blattman2002estimating}\\
  $\theta_C$ & Level of selectivity of clonal competition amongst tumour cells & [0.1, 2] & \textit{ad hoc}\\
  $\theta_T$ & Level of selectivity of self-regulation mechanisms of CTLs & [0.1, 2] & \textit{ad hoc}\\
  $\gamma$ & TCR-tumour antigen binding affinity & [0.1, 3.5] & \citep{stromberg2006robustness,stromberg2013diversity}\\
  $\lambda_C$ & Probability of phenotypic variation of tumour cells & 0.01 & \citep{stace2019phenotype}\\
   %$\chi$ & Phenotype-step & 0.0013 & \textit{ad hoc}\\
   %$\tau$ & Time-step & 0.05 & \textit{ad hoc}\\
   %$t_f$ & Final time & 30 days & \textit{ad hoc}\\
  
  \hline
\end{tabular}
  \label{Tab:Tcr}
\end{table}

\subsection{Main results}
\label{Main results}
\begin{figure}[h!]
\centering 
\includegraphics[scale=0.6]{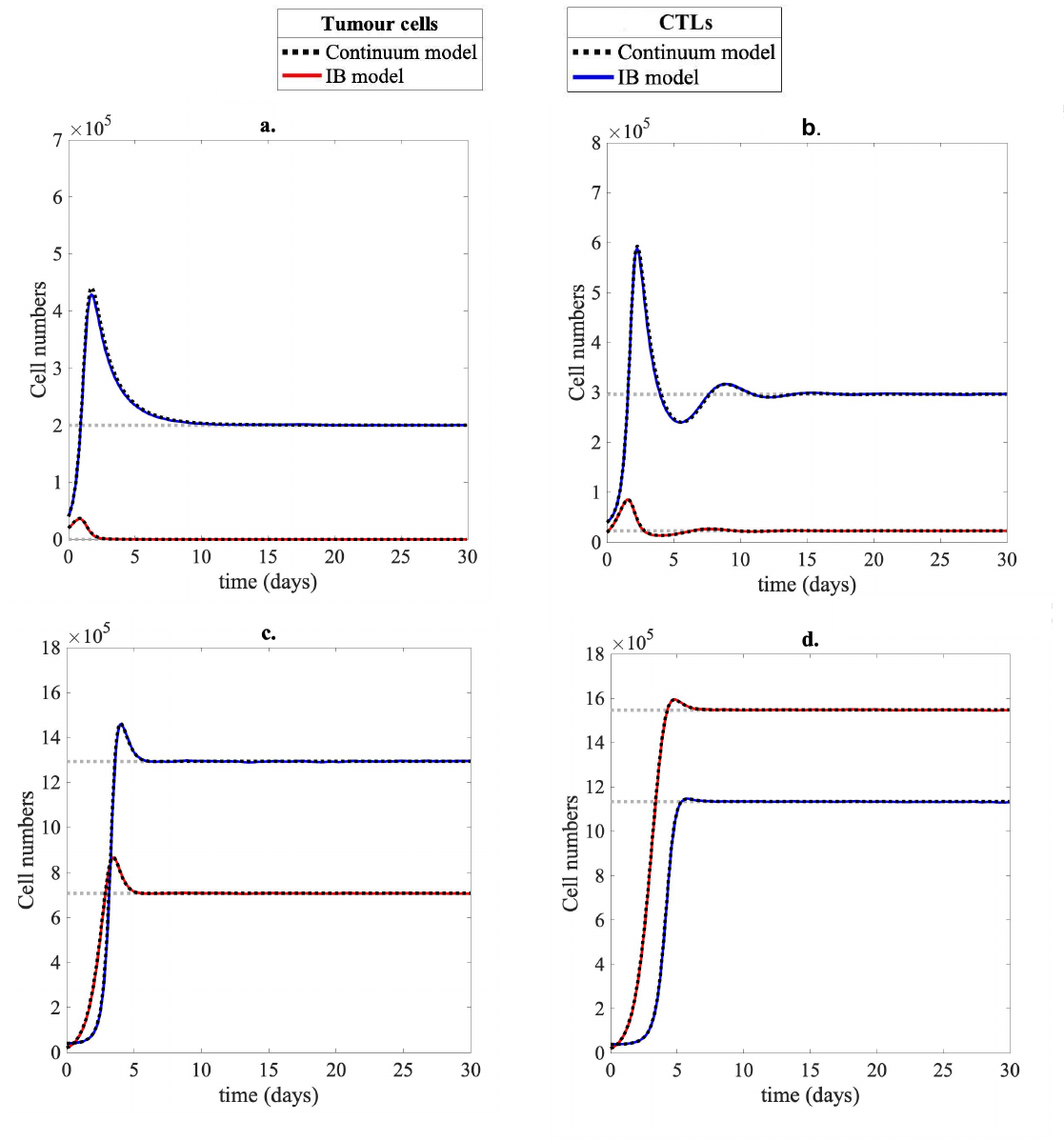}
\caption{\textbf{Eradication of tumour cells and emergence of hot tumour-like, altered tumour-like and cold tumour-like scenarios.} Panel \textbf{a.} displays the plots of the time evolution of the total number of tumour cells ($\rho_C$) and CTLs ($\rho_T$) of the individual-based model (solid, coloured lines) and the continuum model (dashed, black lines) when $\gamma$ is high enough that condition~\eqref{condition1} is satisfied (i.e. condition~\eqref{condition} does not hold). Here, $\alpha_T=0.5$/day and all the other parameters are as in Table~\ref{Tab:Tcr} with $\gamma=3.5$, $\eta=1.8$ and $\theta_C=\theta_T=1.8$. The grey dotted lines highlight the steady-state values of $\rho_C$ and $\rho_T$ given by~\eqref{steady_b}. Panels \textbf{b.}-\textbf{d.} display similar plots for sufficiently large, intermediate and sufficiently small values of $\gamma$ that satisfy condition~\eqref{condition} -- i.e. $\gamma=2$ (panel \textbf{b.}), $\gamma=0.3$ (panel \textbf{c.}) and $\gamma=0.12$ (panel \textbf{d.}). All the other parameters are as in Table~\ref{Tab:Tcr} with $\eta=1.8$ and $\theta_C=\theta_T=1.8$. The grey dotted lines highlight the steady-state values of $\rho_C$ and $\rho_T$ given by~\eqref{steady_d}. Initial conditions~\eqref{initcondAgentbased} and~\eqref{initcondPDE} with $a=0$ were used to carry out numerical simulations. Analogous results hold when $a>0$ in~\eqref{initcondAgentbased} and~\eqref{initcondPDE} (\emph{cf.} Figure~\ref{fig:new1bis} in Appendix~\ref{AppendixS}). The results from the individual-based model correspond to the average over two realisations and the related variance is displayed by the coloured areas surrounding the curves.}
\label{fig:new1}
\end{figure}

\paragraph{Eradication of tumour cells} When $\gamma$ is high enough so that condition~\eqref{condition1} is satisfied (i.e. condition~\eqref{condition} does not hold), after initial growth, the total number of tumour cells decreases steadily over time until the tumour cell population is completely eradicated (\emph{cf.} Figure~\ref{fig:new1}{\bf a}). This is due to the fact that, in response to a rapid growth in the size of the tumour cell population, the high TCR-tumour antigen binding affinity allows the population of CTLs to embark on rapid expansion in size that continues until CTLs have reached the critical mass required to push the population of tumour cells towards extinction. The expansion of the CTL population is followed by the transition to a contraction phase, which is characterised by a decline of the total number of CTLs to a level corresponding to the maintenance of a form of immunological memory. In fact, CTLs can persist after tumour eradication and could develop into memory T cells, thus preventing tumour outgrowth~\citep{prehn1957immunity,zhang2007induction}.  

\paragraph{Hot tumour-like scenarios} When $\gamma$ satisfies condition~\eqref{condition} but is still sufficiently high, the total number of CTLs attains a value large enough to keep the total number of tumour cells steadily low. After initial growth, the total number of tumour cells decreases over time until it stabilises itself around a relatively small value (\emph{cf.} Figure~\ref{fig:new1}{\bf b}). As a result, the average value of the immune score $\overline{I}$ defined via~\eqref{immunoscore} and~\eqref{def:avgI} is one order of magnitude larger than $1$ (i.e. for the parameter values considered here $\overline{I} \approx 12.7$). In the framework of our model, this corresponds to the emergence of \emph{hot tumour-like scenarios}. 

\paragraph{Altered tumour-like scenarios} For intermediate values of $\gamma$ that satisfy condition~\eqref{condition}, after initial growth, a certain number of tumour cells and a slightly larger number of CTLs stably coexist (\emph{cf.} Figure~\ref{fig:new1}{\bf c}). In this case, the average value of the immune score $\overline{I}$ defined via~\eqref{immunoscore} and~\eqref{def:avgI} is just slightly larger than $1$ (i.e. for the parameter values considered here $\overline{I} \approx 1.6$). In the framework of our model, this corresponds to the emergence of \emph{altered tumour-like scenarios}.  

\paragraph{Cold tumour-like scenarios} For sufficiently small values of $\gamma$ that satisfy condition~\eqref{condition}, in the early stage of cell dynamics the total number of tumour cells overtakes the total number of CTLs, and keeps expanding until saturation (\emph{cf.} Figure~\ref{fig:new1}{\bf d}). Accordingly, the average value of the immune score $\overline{I}$ defined via~\eqref{immunoscore} and~\eqref{def:avgI} is smaller than $1$ (i.e. for the parameter values considered here $\overline{I} \approx 0.7$), which corresponds to the emergence of \emph{cold tumour-like scenarios} in the framework of our model. 

\begin{figure}[h!]
\centering 	
\includegraphics[scale=0.7]{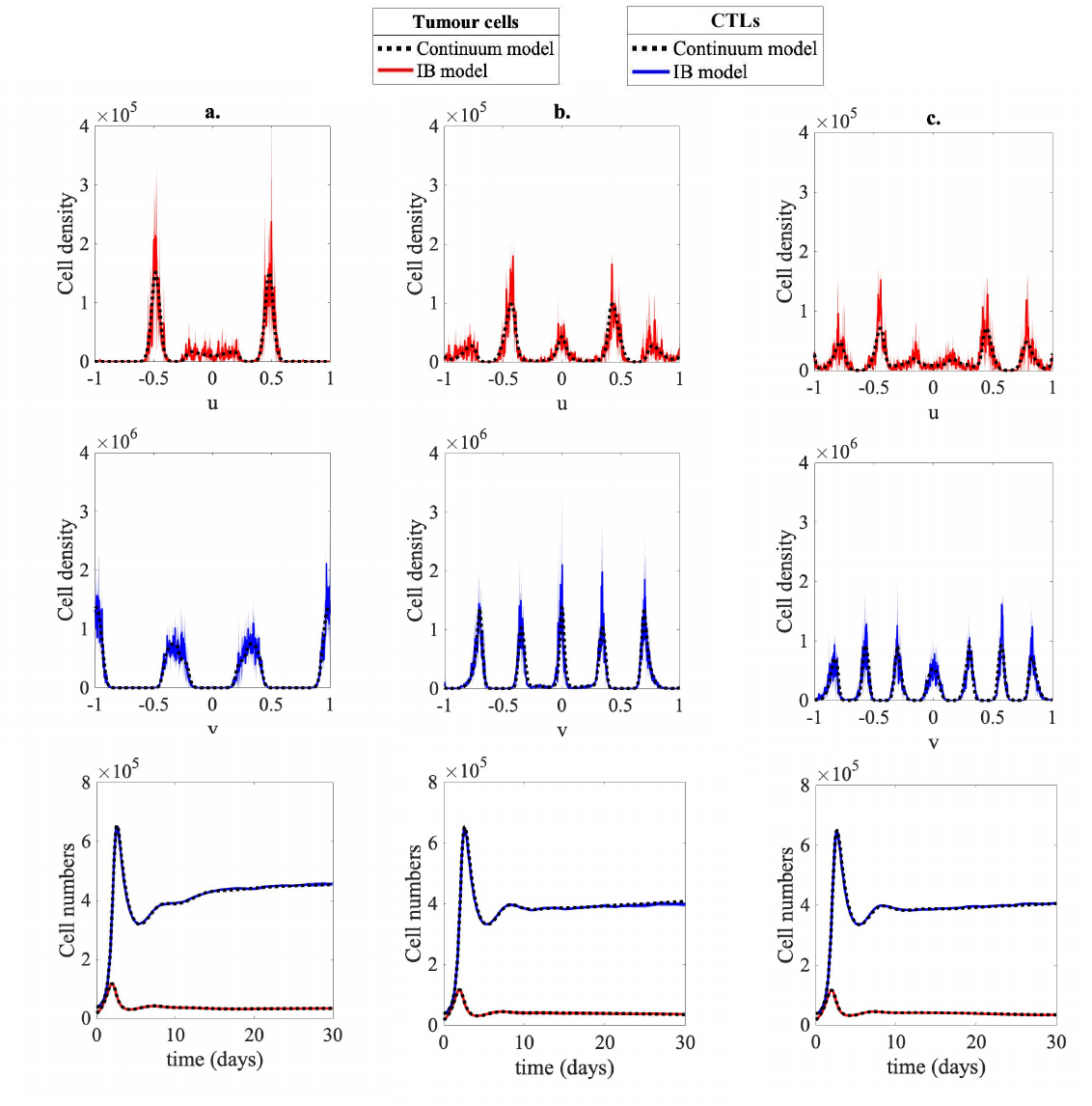}
\caption{\textbf{Patterns of phenotypic coevolution between tumour cells and CTLs: impact of the parameters $\theta_C$ and $\theta_T$.} Top panels display the plots of the population density function of tumour cells ($n_C$) and central panels display the plots of the population density function of CTLs ($n_T$) of the individual-based model (solid, coloured lines) and continuum model (dashed, black lines) at the end of numerical simulations (i.e. at $t=30$) when conditions~\eqref{condition} and~\eqref{cond:patterns} are satisfied and progressively smaller values of $\theta_C$ and $\theta_T$ are considered -- i.e. $\theta_C = \theta_T = 0.5$ (panels \textbf{a.}), $\theta_C = \theta_T = 0.3$ (panels \textbf{b.}) and $\theta_C = \theta_T = 0.2$ (panels \textbf{c.}). All the other parameters are as in Table~\ref{Tab:Tcr} with $\gamma=1.5$ and $\eta=0.7$. Bottom panels display the corresponding plots of the time evolution of the total number of tumour cells ($\rho_C$) and CTLs ($\rho_T$). Initial conditions~\eqref{initcondAgentbased} and~\eqref{initcondPDE} with $a=1$ and $A=5$ were used to carry out numerical simulations. Analogous results were obtained when using different values of the parameter $A$ (results not shown). The results from the individual-based model correspond to the average over two realisations of the underlying  random walk and the related variance is displayed by the coloured areas surrounding the curves.
\label{fig:new2}} 
\end{figure}

\begin{figure}[h!]
\centering 	
\includegraphics[scale=0.7]{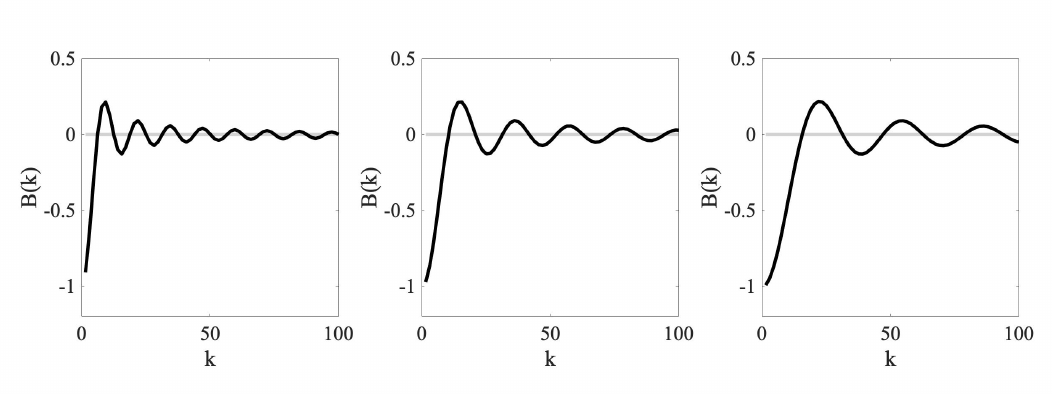}
\caption{\textbf{Plots of $B$ for different values of $\theta_C$ and $\theta_T$.} Plots of the function $B(k,x)$ defined via~\eqref{eq:B(k)} for any $x \in \displaystyle{\argmin_{x \in [-1,1]}} \Psi(x;\theta)$ with $\theta = \theta_C = \theta_T$, under the parameter values used in Figure~\ref{fig:new2} -- i.e. $\theta = \theta_C = \theta_T = 0.5$ (left panel), $\theta = \theta_C = \theta_T = 0.3$ (central panel) and $\theta = \theta_C = \theta_T = 0.2$ (right panel).
\label{fig:final3}} 
\end{figure}

\paragraph{Robustness of numerical results}
The plots in Figure~\ref{fig:new1} demonstrate that there is an excellent quantitative agreement between the results of numerical simulations of the individual-based model and numerical solutions of the corresponding continuum model. Moreover, consistently with the results of linear stability analysis of the continuum model presented in Section~\ref{temporal perturbations}, these numerical results show that the total numbers of tumour cells and CTLs converge either to the steady-state values given by~\eqref{steady_b} (\emph{cf.} Figure~\ref{fig:new1}{\bf a}), or the steady-state values given by~\eqref{steady_d} (\emph{cf.} Figure~\ref{fig:new1}{\bf b-d}), depending on the fact that the choices of the model parameters are such that condition~\eqref{condition1} or condition~\eqref{condition} holds, respectively. When convergence to the steady state $(\rho_{C2}^*,\rho_{T2}^*)$ given by~\eqref{steady_d} occurs, in the long run, the value of the average immune score $\overline{I}$ defined via~\eqref{def:avgI} reflects the value of the ratio $\rho_{T2}^*/\rho_{C2}^*$. Therefore, in the framework of our tumour classification based on the average immune score $\overline{I}$ (see page 5), if condition~\eqref{condition} is met: cold tumour-like scenarios and hot tumour-like scenarios will emerge when the values of the model parameters are such that the ratio $\rho_{T2}^*/\rho_{C2}^*$ is smaller than $1$ or at least about one order of magnitude larger than $1$, respectively, whereas altered tumour-like scenarios will emerge in the remaining cases. This has been confirmed by the results of additional numerical simulations (results not shown). Hence, independently of the specific values of the model parameters, provided that assumption~\eqref{condition} is satisfied, cell dynamics qualitatively similar to those of Figure~\ref{fig:new1}, and corresponding to hot, altered or cold tumour scenarios, will be observed depending on the value of the ratio $\rho_{T2}^*/\rho_{C2}^*$. This testifies to the robustness of the numerical results presented here.
\paragraph{Patterns of phenotypic coevolution between tumour cells and CTLs: impact of the parameters $\theta_C$ and $\theta_T$} Figure~\ref{fig:new2} displays the plots of the phenotype distributions of tumour cells (top panel) and CTLs (central panel) at the end of numerical simulations (i.e. close to numerical equilibrium) alongside the plots of the corresponding time evolution of the total cell numbers (bottom panel). In agreement with the analytical results presented in Section~\ref{spatio-temporal perturbations}, when condition~\eqref{condition} is satisfied and conditions~ \eqref{cond:thetaeta2} and~\eqref{cond:patterns} are met as well, patterns of phenotypic coevolution between tumour cells and CTLs may emerge. Moreover, the top and central panels of Figure~\ref{fig:new2} show that, coherently with the shape of the function $B$ defined via~\eqref{eq:B(k)} (\emph{cf.} the plots in Figure~\ref{fig:final3}), smaller values of $\theta_C$ and $\theta_T$ correlate with the formation of more peaks in the phenotype distributions of the two cell populations. The plots in Figure~\ref{fig:new2} also demonstrate that there is a good agreement between numerical simulations of the individual-based and continuum models. 

\begin{figure}[h!]
\centering 	
\includegraphics[scale=0.4]{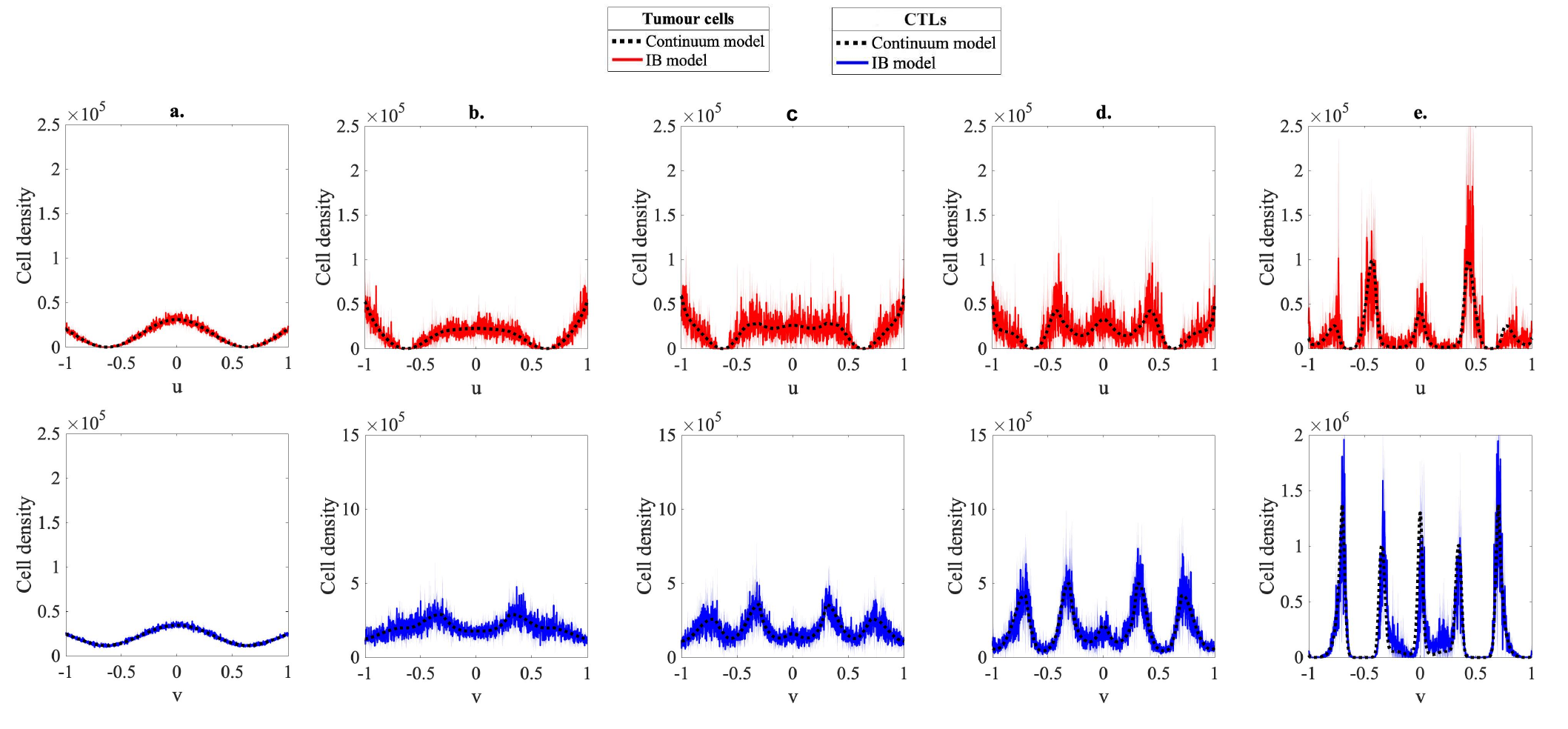}
\caption{\textbf{Sample temporal dynamics of patterns of phenotypic coevolution between tumour cells and CTLs.} Top panels display the plots of the population density function of tumour cells ($n_C$) and bottom panels display the plots of the population density function of CTLs ($n_T$) of the individual-based model (solid, coloured lines) and continuum model (dashed, black lines) at five successive time instants -- i.e. $t=0.4$ (panels \textbf{a.}), $t=4$ (panels \textbf{b.}), $t=10$ (panels \textbf{c.}), $t=16$ (panels \textbf{d.}), $t=30$ (panels \textbf{e.}) -- in the case where condition~\eqref{condition} is satisfied. Here, $\theta_C=\theta_T=0.3$, $\gamma=1.5$ and $\eta=0.7$, and all the other parameters are as in Table~\ref{Tab:Tcr}. Initial conditions~\eqref{initcondAgentbased} and~\eqref{initcondPDE} with $a=1$ and $A=5$ were used to carry out numerical simulations. Analogous results were obtained when using different values of the parameter $A$ (results not shown). The results from the individual-based model correspond to the average over two realisations and the related variance is displayed by the coloured areas surrounding the curves.\label{fig:final4}}
\end{figure}

Sample temporal dynamics of such patterns are summarised by the plots in Figure~\ref{fig:final4}, which show that clonal expansion leads to a rapid proliferation of CTLs that are targeted to the antigens mostly expressed by tumour cells, whereas self-regulation mechanisms induce formerly stimulated CTLs to decay. In turn, the antigen-specific cytotoxic action of CTLs causes the selection of those tumour cells that are able to escape immune recognition. As a result, immune competition induces the formation of multiple peaks in the phenotype distribution of tumour cells. This concurrently shapes the phenotype distribution of CTLs with a time shift corresponding to the time required for the CTLs to adapt to the antigenic distribution of tumour cells. The plots in Figure~\ref{fig:final4} demonstrate that there is again a good agreement between numerical simulations of the individual-based and continuum models. 

\begin{figure}[h!]
\centering 	
\includegraphics[scale=0.5]{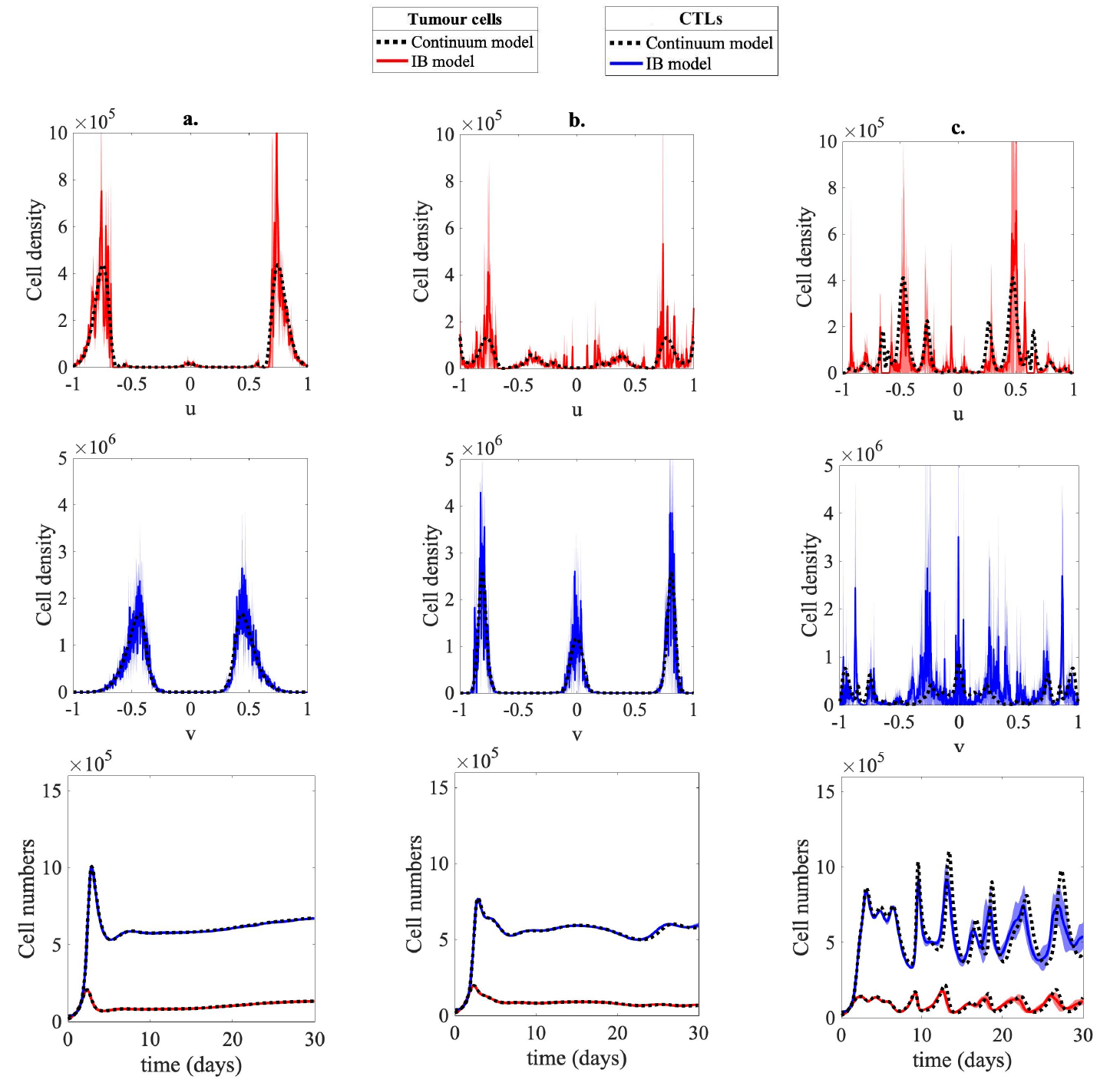}
\caption{\textbf{Patterns of phenotypic coevolution between tumour cells and CTLs: impact of the parameter $\eta$.} Top panels display the plots of the population density function of tumour cells ($n_C$) and bottom panels display the plots of the population density function of CTLs ($n_T$) of the individual-based model (solid, coloured lines) and continuum model (dashed, black lines) at the end of numerical simulations (i.e. at $t=30$) when condition~\eqref{condition} is satisfied and progressively smaller values of $\eta$ are considered -- i.e. $\eta=1$ (panels \textbf{a.}), $\eta=0.6$ (panels \textbf{b.}) and $\eta=0.2$ (panels \textbf{c.}). All the other parameters are as in Table~\ref{Tab:Tcr} with $\gamma=1$ and $\theta_C=\theta_T=0.7$. Bottom panels display the corresponding plots of the time evolution of the total number of tumour cells ($\rho_C$) and CTLs ($\rho_T$). Initial conditions~\eqref{initcondAgentbased} and~\eqref{initcondPDE} with $a=1$ and $A=5$ were used to carry out numerical simulations. Analogous results were obtained when using different values of the parameter $A$ (results not shown). The results from the individual-based model correspond to the average over five realisations of the underlying  random walk and the related variance is displayed by the coloured areas surrounding the curves. \label{fig:final5}} 
\end{figure}

%Under the parameter choice of the simulations of Figure~\ref{fig:final5}, we have the emergence of altered tumour-like scenarios (\emph{cf.} bottom panels of Figure~\ref{fig:final5}), with $\overline{I} \approx 2.8$ in Figure~\ref{fig:final5}a  and Figure~\ref{fig:final5}b and $\overline{I}\approx 2.9$ in Figure~\ref{fig:final5}c. 
\paragraph{Patterns of phenotypic coevolution between tumour cells and CTLs: impact of the parameter $\eta$} The results of numerical simulations summarised by the plots in Figure~\ref{fig:final5} extend the analytical results presented in Section~\ref{spatio-temporal perturbations} by showing that, when condition~\eqref{condition} is satisfied and $\eta$ is sufficiently small, smaller values of $\eta$ may induce the formation of patterns of phenotypic coevolution between tumour cells and CTLs whereby less regular multi-peaked phenotype distributions of the two cell populations emerge (\emph{cf.} top and central panels of Figure~\ref{fig:final5}). The temporal dynamics of such patterns are qualitatively similar to those presented in Figure~\ref{fig:final4} (results not shown). Moreover, numerical simulations indicate that smaller values of $\eta$ correlate with the emergence of oscillations in the total numbers of tumour cells and CTLs, that is, CTLs undergo a succession of expansion and contraction phases that result in an alternate decay and growth of tumour cells~(\emph{cf.} bottom panel in Figure~\ref{fig:final5}c). The plots in Figure~\ref{fig:final5} demonstrate that there is a good agreement between numerical simulations of the individual-based and continuum models. 

\begin{figure}[h!]
\centering 	
\includegraphics[scale=0.62]{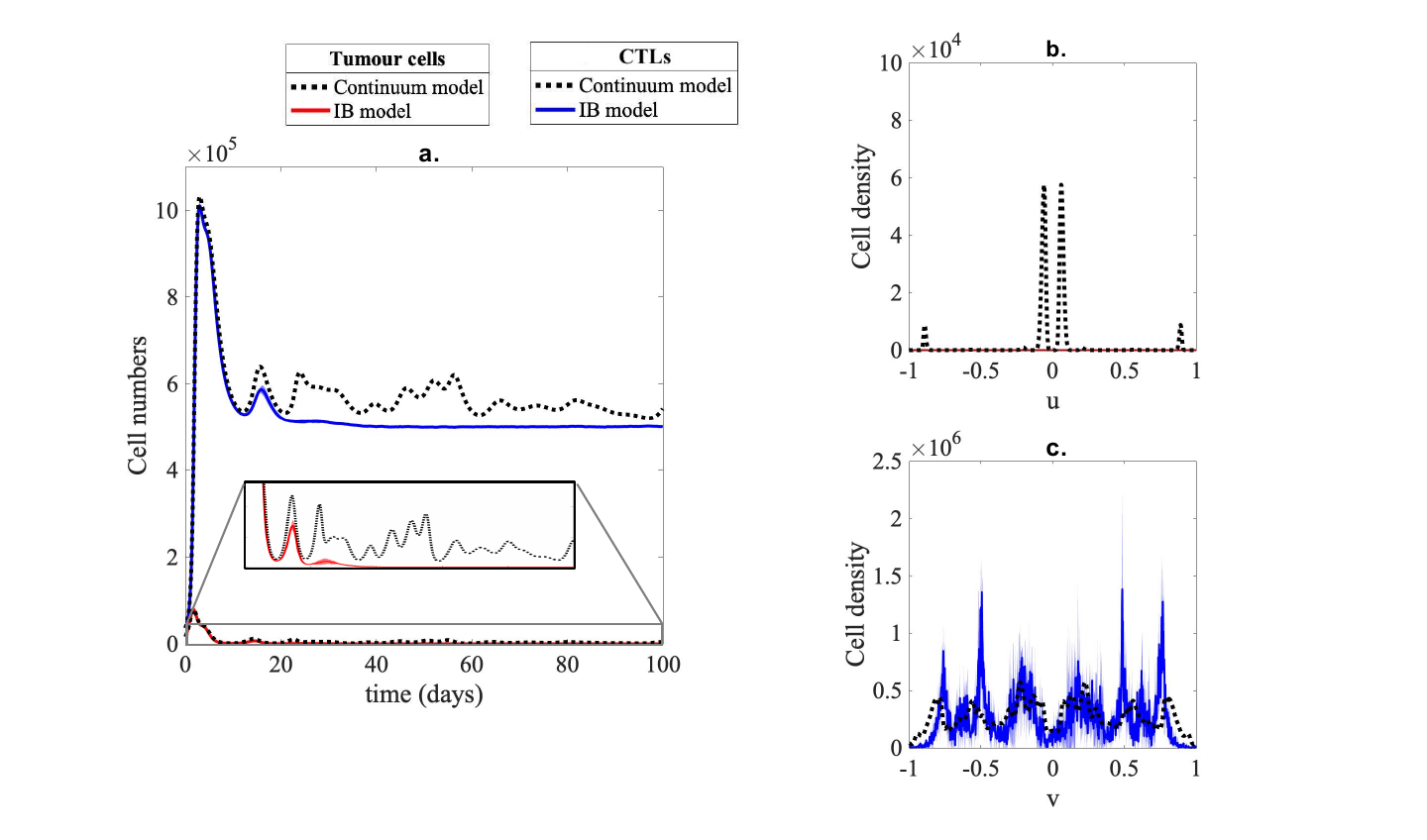}
\caption{\textbf{Possible discrepancies between individual-based and continuum models.} Panel \textbf{a.} displays the plot of the time evolution of the total number of tumour cells ($\rho_C$) and CTLs ($\rho_T$) of the individual-based model (solid, coloured lines) and the continuum model (dashed, black lines) when condition~\eqref{condition} holds, evolution towards relatively small tumour cell numbers occurs, and the parameter $\eta$ is sufficiently small so that less regular multi-peaked cell phenotype distributions emerge -- i.e. $\alpha_T=0.5$/day, $\mu_T=2 \times 10^{-6}\mu l$/day, $\gamma=1.1$, $\eta=0.1$, and all the other parameters as in Table~\ref{Tab:Tcr} with $\theta_C=\theta_T=1.8$. The plots in panels \textbf{b.} and \textbf{c.} display the corresponding population density functions of tumour cells ($n_C$) and CTLs ($n_T$) of the individual-based model (solid, coloured lines) and of the continuum model (dashed, black lines) at the end of simulations (i.e. at $t=t_f=100$). Initial conditions~\eqref{initcondAgentbased} and~\eqref{initcondPDE} with $a=1$ and $A=5$ were used to carry out numerical simulations. Analogous results were obtained when using different values of the parameter $A$ (results not shown). The results from the individual-based model correspond to the average over five realisations and the related variance is displayed by the coloured areas surrounding the curves.} 
\label{fig:new6}
\end{figure}

\paragraph{Possible discrepancies between individual-based and continuum models} The results that have been presented so far indicate that there is a good agreement between the results of computational simulations of the individual-based model and the numerical solutions of the corresponding continuum model. However, we expect possible differences between the outputs of the two models to emerge in the presence of lower tumour cell numbers, which may lead to more pronounced demographic stochasticity, and less regular multi-peaked cell phenotype distributions, which may cause a reduction in the quality of the approximations employed in the formal derivation of the deterministic continuum model from the individual-based model. In order to investigate this, we carried out numerical simulations of the two models for choices of parameter values such that condition~\eqref{condition} holds, evolution towards relatively small tumour cell numbers occurs, and less regular cell phenotype distributions with multiple peaks emerge (see caption of Figure~\ref{fig:new6} for more details). The results obtained are summarised by the plots in Figure~\ref{fig:new6}, which show that the individual-based model predicts eradication of the tumour cell population, whereas the continuum model predicts coexistence between tumour cells and CTLs.
 \section{Discussion, conclusions and research perspectives} 	 
 \label{conclusion}
 \paragraph{Discussion and conclusions} We developed an individual-based model for the coevolutionary  dynamics  between tumour cells and CD8+ cytotoxic T lymphocytes that takes into account the selectivity of antigen-specific immunity. We formally derived the deterministic continuum counterpart of such an individual-based model, and we integrated the results of numerical simulations of the two models with the results of steady-state and linear-stability analyses of the continuum model equations.
 
The results presented in this study shed light on the way in which different parameters shape the coevolutionary dynamics between tumour cells and CD8+ cytotoxic T lymphocytes. In particular, we demonstrated that, \emph{ceteris paribus}, higher values of the TCR-tumour antigen binding affinity (i.e. the parameter $\gamma$ in the model) promote the eradication of tumour cells by CTLs, while lower values facilitate the coexistence between tumour cells and CTLs. Specifically, progressively reducing the TCR-tumour antigen binding affinity brings about the emergence of: \emph{hot tumour-like scenarios}, which are characterised by a large number of \textit{in situ} CTLs and a low number of tumour cells, and thus represent a more fertile ground for anticancer therapeutic intervention; \emph{altered tumour-like scenarios}, which reflect the intrinsic ability of the immune system to effectively mount a CTL-mediated immune response and the ability of tumour cells to partially escape such a response; \emph{cold tumour-like scenarios}, which are characterised by an insufficient number of \textit{in situ} CTLs and are invariably associated with poor prognosis~\citep{galon2019approaches}. This
classification of tumours is also supported by experimental works showing that \textit{in situ} immune reaction might be the strongest parameter influencing clinical outcome, regardless of the local tumour extension and its spread to lymph nodes~\citep{pages2008essential,galon2006type, galon2016immunoscore}. Moreover, our findings support the idea that TCR-tumour antigen binding affinity may be a good intervention target for immunotherapy that aims to turn cold or altered tumours into hot ones by enhancing CTL response. In this regard, our findings are in agreement with the conclusions of previous experimental articles indicating that a strong binding affinity of T cells to tumour antigens may play a key role in the overall immune response to the disease \citep{gerdemann2011cytotoxic}.In particular, in altered tumours, increasing antigenicity, via the removal of co-inhibitory signals and/or the supply of co-stimulatory signals~\citep{galon2016immunoscore, whiteside2016emerging}, may enhance \textit{in situ} CTLs activity, and has proven to be effective in the treatment of advanced-stage melanoma~\citep{wolchok2017overall}, renal cell carcinoma~\citep{motzer2018nivolumab} and non-small cell lung cancer~\citep{hellmann2018nivolumab}. In cold tumours, a proposed
approach to overcome the lack of a pre-existing immune response consists in combining
a priming therapy that boosts CTL responses with
the removal of co-inhibitory signals through approaches such as immune checkpoint~\citep{galon2019approaches}. The therapeutic success achieved by combining immune checkpoint therapy with chemotherapy in metastatic NSCLC has demonstrated the potential strength of this dual approach~\citep{gandhi2018pembrolizumabe}.

Moreover, the results presented here indicate that the affinity range of TCRs (i.e. the parameter $\eta$ in the model), the selectivity of clonal competition amongst tumour cells (i.e. the inverse of the parameter $\theta_C$ in the model) and the selectivity of self-regulation mechanisms acting on CD8+ cytotoxic T lymphocytes (i.e. the inverse of the parameter $\theta_T$ in the model) play a pivotal role in the formation of patterns of phenotypic coevolution, which create the substrate for the emergence of less regular cell phenotype distributions with multiple peaks. Such patterns are underpinned by some form of immunoediting whereby the population of CTLs evolves and continuously adapts its receptor repertoire in order to recognise and effectively eliminate tumour cells and, in turn, the antigen-specific selective pressure exerted by CTLs leads to the selection of those tumour clones that are able to evade immune recognition~\citep{dunn2002cancer}. The adaptability of tumour cells and CTLs and the selective pressure they mutually exert on each other during cancer development are emerging as crucial factors in determining cancer evolutionary trajectories. This has been shown in the context of chronic lymphocytic leukemia \citep{purroy2017coevolution} and other cancer types, as reviewed in~\citep{george2021implications}. Our results offer also a theoretical basis for the development of anti-cancer therapy aiming at engineering TCRs so as to shape their affinity for cancer targets~\citep{border2019affinity,crean2020molecular,li2019genetically,zhao2021engineered} and adaptive therapy aiming at altering intratumour clonal competition~\citep{gatenby2009adaptive,west2020towards}, in order to control the coevolutionary dynamics between tumour cells and CD8+ cytotoxic T lymphocytes. In this respect, one of the
best known treatment based on engineering specific TCRs is based on CAR-T
cells \citep{ wang2016clinical}, which confer CTLs the ability to
target specific antigens. It has been demonstrated that this therapeutic strategy has several potential advantages over conventional therapies, including specificity, rapidity, high success rate and long-lasting effects~\citep{june2018car,gomes2018cancer}.

The good agreement between the results of numerical simulations of the individual-based and continuum models, along with the quantitative information given by~\eqref{steady_b} and~\eqref{steady_d} and the precise conditions given by~\eqref{condition1} and~\eqref{cond:patterns}, testifies to the robustness of the  biological insight gained in this work. We also showed that possible differences between cell dynamics produced by the individual-based and continuum models can emerge under parameter settings that correspond to less regular cell phenotype distributions and more pronounced demographic stochasticity. In fact, these cause a reduction in the quality of the approximations employed in the formal derivation of the deterministic continuum model from the individual-based model (cf. Appendix \ref{Appendix A}). This demonstrates the importance of integrating individual-based and continuum approaches when considering mathematical models for tumour-immune competition.

\paragraph{Research perspectives} From a mathematical point of view, we plan to carry out a systematic investigation of the conditions on the affinity range of TCRs that may lead to the emergence of
oscillations in cell numbers observed in the numerical simulations presented in this work. Moreover, from a modelling point of view, our individual-based modelling framework for the coevolutionary dynamics  between tumour cells and CD8+ cytotoxic T lymphocytes, along with the formal derivation of the corresponding continuum  model,  can be developed further in several ways. For instance, a myriad of immunosuppressive strategies, the so-called immune checkpoints, help tumour cells acquiring features that enable them to evade immune detection, which may ultimately induce the exhaustion of CTLs in the tumour micro-environment, which impairs the immune response. The modelling approach presented here does not capture this aspect. However, exhaustion mechanisms could be incorporated into the individual-based model by, for example, allowing CTLs to enter a suppressed state (i.e. CTLs would become exhausted and thus would no longer able to eliminate tumour cells). In the continuum model, this would result in the presence of an additional loss term in the IDE~\eqref{eq:PDEs}$_2$ along with a third equation for the dynamics of exhausted CTLs. Another track to follow to further enrich our model would be to include a spatial structure, for instance by embedding the tumour cells in the geometry of a solid tumour, and to take explicitly into account the effect of both spatial and antigen-specific interactions between tumour cells and CTLs, as similarly done in~\citep{kather2017silico, macfarlane2018modelling, macfarlane2019stochastic}. Including a spatial structure would make it possible, \emph{inter alia}, to introduce a more precise definition of the immune score that incorporates the level of CTL infiltration. Furthermore, at this stage, the mathematical representation of the phenotypic state of tumour cells and CTLs employed in our modelling framework is rather abstract. This might make it difficult to carry out precise quantitative comparisons between the results of numerical simulations and experimental data. This limitation could be overcome by employing a mathematical representation of tumour antigens and TCRs similar to the one that we proposed in~\citep{leschiera2021mathematical}, whereby a discrete set of tumour antigens that can be recognised by a unique repertoire of TCRs is considered. Finally, it would be interesting to incorporate explicitly into the model the effects of immunotherapeutic agents or other therapeutic agents. These are all lines of research that we will be pursuing in the future.  

\section*{Funding} 
E.L. has received funding from the European Research Council (ERC) under the European Union's Horizon2020 research and innovation programme (grant agreement No 740623). T.L. gratefully acknowledges support from the Italian Ministry of University and Research (MUR) through the grant ``Dipartimenti di Eccellenza 2018-2022'' (Project no. E11G18000350001) and the PRIN 2020 project (No. 2020JLWP23) ``Integrated Mathematical Approaches to Socio–Epidemiological Dynamics'' (CUP: E15F21005420006). L.A., E.L. and T.L. gratefully acknowledge support from the CNRS International Research Project ``Modélisation de la biomécanique cellulaire et tissulaire'' (MOCETIBI).

\newpage
\begin{appendices}
\section{Formal derivation of the continuum model}
  	   			\label{Appendix A}
Using a method analogous to that employed in~\citep{ardavseva2020comparative, chisholm2016evolutionary, stace2019phenotype}, we show that the PDE-IDE system~\eqref{eq:PDEs} can be formally derived as the appropriate continuum limit of the individual-based model presented in this article. 

Substituting definitions~\eqref{pC} of ${\rm P}^b_C$ and ${\rm P}^q_C$ into the difference equation~\eqref{eq:mastereqs}$_1$ for $n_{C_{i}}^{h}$ and definitions~\eqref{pT} of ${\rm P}^b_T$ and ${\rm P}^q_T$ into the difference equation~\eqref{eq:mastereqs}$_2$ for $n_{T_{j}}^{h}$ yields
\begin{equation}
\label{eq:mastereqsrev}
\begin{cases}
n_{C_{i}}^{h+1}=\left[1 + \tau \, \alpha_C - \tau \, \left(\mu_C \, K_{C_{i}}^{h} + \zeta_C \,\gamma\, J_{C_{i}}^{h}\right) \right]\left[\dfrac{\lambda_C}{2} \left(n_{C_{i+1}}^{h} + n_{C_{i-1}}^{h}\right) + (1-\lambda_C) \, n_{C_{i}}^{h}\right],
\\\\
n_{T_{j}}^{h+1}=\left[1 + \tau \, \left(\alpha_T + \zeta_T \,\gamma\, J_{T_{j}}^{h}\right) -  \tau \, \mu_T \, K_{T_{j}}^{h} \right] \, n_{T_{i}}^{h},
\end{cases}
\end{equation} 
where $n_{C_{i}}^{h} \equiv n_{C}(u_i, t_h)$ with $(u_i, t_h) \in \mathcal{I} \times (0,t_f]$ and $n_{T_{j}}^{h} \equiv n_{T}(v_j, t_h)$ with $(v_j, t_h) \in \overline{\mathcal{I}} \times (0,t_f]$. Using the fact that the following relations hold for $\tau$ and $\chi$ sufficiently small
	    		$$t_h\approx t,\quad t_{h+1}\approx t+\tau,\quad u_i\approx u,\quad u_{i\pm 1}\approx u\pm\chi,\quad 
	    		 v_j\approx v,$$
	    		 $$n_{C_{i}}^{h}\approx n_C(u,t), \quad n_{C_{i}}^{h+1}\approx n_C(u,t+\tau), \quad n_{C_{i\pm 1}}^{h}\approx n_C(u\pm \chi,t), \quad \rho_C^h \approx \rho_C(t):=\int_{\mathcal{I}}n_C(u,t) \, \mathrm{d}u,$$
	    		 $$ J_{C_{i}}^{h}\approx J_C(u,t):=\int_{\mathcal{I}} g(u,v; \eta) \, n_T(v,t) \, \mathrm{d}v, \quad K_{C_{i}}^{h}\approx K_C(u,t):=\int_{\mathcal{I}} g(u,w; \theta_C) \, n_C(w,t) \, \mathrm{d}w,$$
	    		 
	    		 $$n_{T_{j}}^{h}\approx n_T(v,t), \quad n_{T_{j}}^{h+1}\approx n_T(v,t+\tau), \quad \rho_T^h \approx \rho_T(t):=\int_{\mathcal{I}}n_T(v,t) \, \mathrm{d}v,$$
	    		 $$ J_{T_{j}}^{h}\approx J_T(v,t):=\int_{\mathcal{I}} g(v,u; \eta) \, n_C(u,t) \, \mathrm{d}u, \quad K_{T_{j}}^{h}\approx K_T(v,t):=\int_{\mathcal{I}} g(v,w; \theta_T) \, n_T(w,t) \, \mathrm{d}w,$$
where the function $g$ is defined via~\eqref{g}, the system of equations~\eqref{eq:mastereqsrev} can be formally rewritten in the approximate form
\begin{equation}
\label{e:deriv1}
\begin{cases}
n_C(u,t+\tau) =\Big[1 + \tau \, R_C(K_C(u,t),J_C(u,t)) \Big] \times 
\\
\qquad \qquad \qquad \qquad \qquad \qquad \times \, 
\left[\dfrac{\lambda_C}{2} \left(n_C(u+ \chi,t) + n_C(u- \chi,t)\right) + (1-\lambda_C) \, n_C(u,t)\right], 
\\\\
n_T(v,t+\tau) =\Big[1 + \tau \, R_T(K_T(v,t),J_T(v,t))\Big] \, n_T(v,t), 
\end{cases}
\end{equation} 
where $u \in \mathcal{I}$, $v \in \overline{\mathcal{I}}$ and $t \in (0,t_f]$. Here,
\begin{equation}
\label{e:defRCRT}
R_C(K_C,J_C):= \alpha_C - \Big(\mu_C \, K_C + \zeta_C \,\gamma\, J_C\Big), \quad R_T(K_T,J_T):= \Big(\alpha_T + \zeta_T \,\gamma\, J_T\Big) -  \mu_T \, K_T. 
\end{equation} 
If the function $n_C(u,t)$ is twice continuously differentiable with respect to the variable $u$, for $\chi$ sufficiently small we can use the Taylor expansions
\begin{equation}\label{taylorexp}
n_C(u\pm \chi,t) = n_C(u,t) \pm \chi \partial_u n_C(u,t) + \frac{\chi^2}{2} \partial^2_{uu} n_C(u,t) + h.o.t. \ .
\end{equation}
Substituting the Taylor expansions~\eqref{taylorexp} into equation~\eqref{e:deriv1}$_1$ for $n_C(u,t+\tau)$, after a little algebra we find
$$
\begin{cases}
\dfrac{n_C(u,t+\tau) - n_C(u,t)}{\tau} - \dfrac{\lambda_C \chi^2}{2 \tau} \partial^2_{uu} n_C(u,t) = R_C(K_C(u,t),J_C(u,t)) \, n_C(u,t) + 
\\
\qquad \qquad \qquad \qquad \qquad \qquad \qquad \qquad \qquad \qquad + \, \dfrac{\lambda_C \chi^2}{2} \, R_C(K_C(u,t),J_C(u,t)) \, \partial^2_{uu} n_C(u,t) + h.o.t. \, , 
\\\\
\dfrac{n_T(v,t+\tau) - n_T(v,t)}{\tau} = R_T(K_T(v,t),J_T(v,t)) \, n_T(v,t).
\end{cases}
$$
If, in addition, the functions $n_C(u,t)$ and $n_T(v,t)$ are continuously differentiable with respect to the variable $t$, letting $\tau \to 0$ and $\chi \to 0$ in such a way that condition~\eqref{conditionBeta} is met, from the latter system of equations we formally obtain
$$
\begin{cases}
\partial_t n_C(u,t) -\beta_C \partial^2_{uu} n_C(u,t) = R_C(K_C,J_C) \, n_C(u,t), \quad (u, t) \in \mathcal{I} \times (0,t_f],
\\\\
\partial_t n_T(v,t) = R_T(K_T,J_T) \, n_T(v,t), \quad (v, t) \in \overline{\mathcal{I}} \times (0,t_f].
\end{cases}
$$
Substituting definitions~\eqref{e:defRCRT} of $R_C(K_C,J_C)$ and $R_T(K_T,J_T)$ into the above system of equations gives the PDE-IDE system~\eqref{eq:PDEs}. Finally, the no-flux boundary conditions~\eqref{neumann} follow from the fact that the attempted phenotypic variation of a tumour cell is aborted if it requires moving into a phenotypic state that does not belong to the interval $\overline{\mathcal{I}}$.
 
\section{Details of numerical simulations of the continuum model}
 \label{Appendix3}
To construct numerical solutions of the PDE-IDE system~\eqref{eq:PDEs} subject both to the no-flux boundary conditions~\eqref{neumann} and to the initial condition~~\eqref{initcondPDE}, we use a uniform discretisation of step $\Delta x = 0.0013$ of the interval $\overline{\mathcal{I}} = [-L,L]$ as the computational domain of the independent variables $u$ and $v$, and a uniform discretisation of step $\Delta t = 0.05$ of the time interval $(0, t_f]$.

We construct numerical solutions of the non-local PDE~\eqref{eq:PDEs}$_1$ for $n_C$ using a time-splitting approach, which is based on the idea of decomposing the original problem into simpler subproblems that are then sequentially solved at each time-step using an explicit Euler method with step $\Delta t$. This leads to the following time-dicretisation of the PDE-IDE system~\eqref{eq:PDEs} subject to the Neumann boundary conditions~\eqref{neumann}:
\begin{equation}
\label{PDE:timedisc}
\begin{cases}
n_C^{k+\frac{1}{2}}(u)=n_C^k(u) + \Delta t \, R_C(K_C^k(u),J_C^k(u)) \, n_C^k(u), \quad u \in [-L,L],
\\\\ 
n_C^{k+1}(u)=n_C^{k+\frac{1}{2}}(u)+\Delta t \, \beta_C \, \partial^2_{uu} n_C^{k+\frac{1}{2}}(u), \quad u \in (-L,L),
\\\\
\partial_u n_C^{k+1}(u) = 0, \quad u \in \{-L,L\}
\\\\
n_T^{k+1}(v)=n_T^k(v) + \Delta t \, R_T(K_T^k(v),J_T^k(v)) \, n_T^k(v), \quad v \in [-L,L],
\end{cases}
\end{equation}
where
$$
R_C(K_C^k,J_C^k):=\alpha_C - \mu_CK_C^k - \zeta_C \gamma J_C^k, \quad R_T(K_T^k,J_T^k):=\alpha_T - \mu_T K_T^k - \zeta_T  \gamma J_T^k.
$$
The system of equations~\eqref{PDE:timedisc} is numerically solved using a three-point finite difference explicit scheme for the diffusion term~\citep{leveque2007finite} and an implicit-explicit finite difference scheme for the remaining terms~\citep{lorenzi2015dissecting, lorz2013populational}, which leads to the following system of equations
$$
\begin{cases}
n_{Ci}^{k+\frac{1}{2}}= n_{Ci}^{k}\dfrac{1+\Delta t \,R_{C}(K_{Ci}^k,J_{Ci}^k)_+}{1+\Delta t\, R_{C}(K_{Ci}^k,J_{Ci}^k)_-}, \quad u_i \in [-L,L],
\\\\ 
n_{Ci}^{k+1} = n_{Ci}^{k+\frac{1}{2}}+\beta_C\Delta t \dfrac{n_{Ci+1}^{k+\frac{1}{2}}-2n_{Ci}^{k+\frac{1}{2}}+n_{Ci-1}^{k+\frac{1}{2}}}{\Delta x^2}, \quad u_i \in (-L,L),
\\\\
n_{C i}^{k+1} = n_{C i-1}^{k+1}, \quad u_i \in \{-L, L\},
\\\\
n_{Tj}^{k+1} = n_{Tj}^{k}\dfrac{1+\Delta t \,R_{T}(K_{Tj}^k,J_{Tj}^k)_+}{1+\Delta t\, R_{T}(K_{Tj}^k,J_{Tj}^k)_-}, \quad v_j \in [-L,L].
\end{cases}
$$
Here, $R_{C}(\cdot,\cdot)_+$ and $R_{T}(\cdot,\cdot)_+$ are the positive parts of $R_C(\cdot,\cdot)$ and $R_T(\cdot,\cdot)$, while $R_{C}(\cdot,\cdot)_-$ and $R_{T}(\cdot,\cdot)_-$ are the negative parts of $R_C(\cdot,\cdot)$ and $R_T(\cdot,\cdot)$. Moreover,
$$
K_{C_{i}}^{k} = \sum_q g(u_i,u_q; \theta_C) \, n_{Cq}^k \, \Delta x, \qquad K_{T_{j}}^{h} = \sum_q g(v_j,v_q; \theta_T) \, n_{Tq}^k \, \Delta x
$$
and
$$
J_{C_{i}}^{h} = \sum_j g(u_i,v_j; \eta) \, n_{Tj}^k \, \Delta x, \qquad J_{T_{j}}^{h} = \sum_i g(v_j,u_i; \eta) \,  n_{Ci}^k \, \Delta x.
$$
Given the values of the parameter $\tau$, $\chi$ and $\lambda_C$ of the individual-based model, the value of the parameter $\beta_C$ is defined so that condition~\eqref{conditionBeta} is met. The other parameter values  are chosen to be coherent with those used to carry out numerical simulations of the individual-based model, which are specified in the main body of the paper.

\newpage
\section{Supplementary figures}
 \label{AppendixS}
 \beginsupplement
 
 \begin{figure}[h!]
\centering 
\includegraphics[scale=0.6]{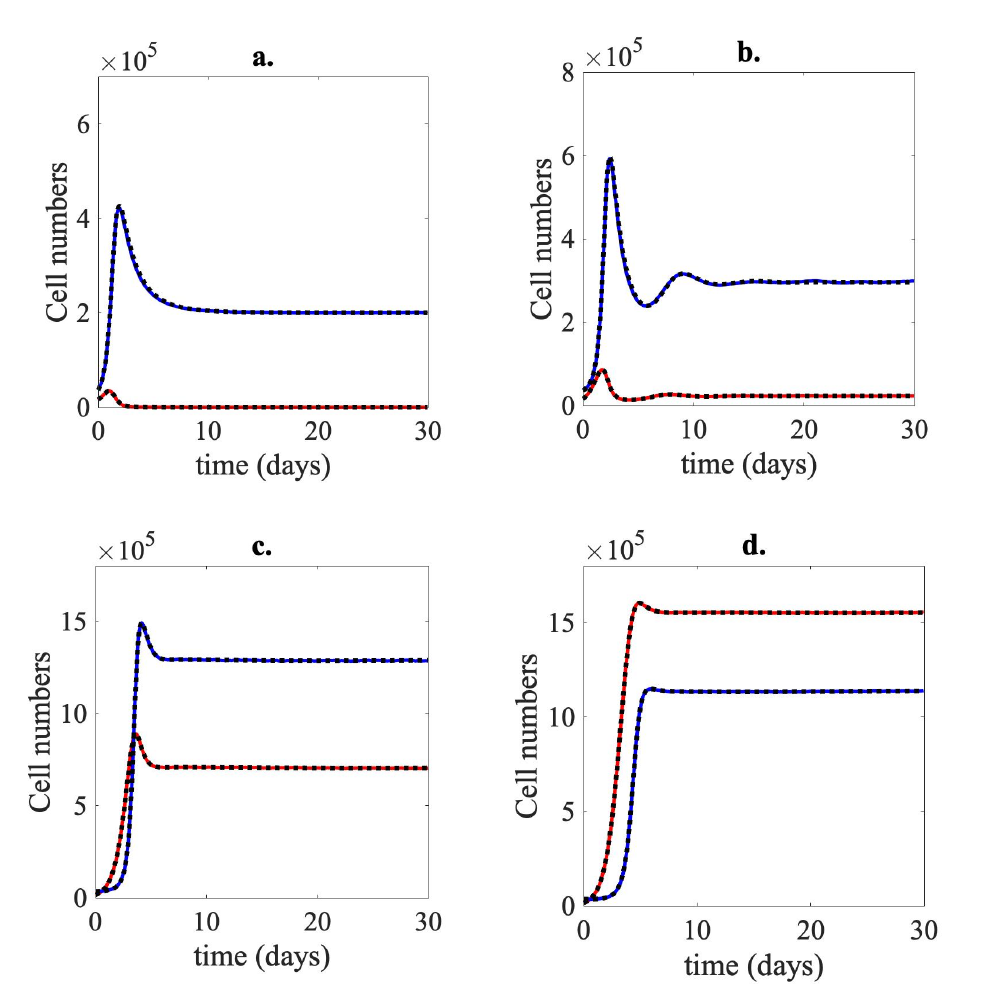}
\caption{\textbf{Eradication of tumour cells and emergence of hot tumour-like, altered tumour-like and cold tumour-like scenarios: the case where $a>0$.} Panel \textbf{a.} displays the plots of the time evolution of the total number of tumour cells ($\rho_C$) and CTLs ($\rho_T$) of the individual-based model (solid, coloured lines) and the continuum model (dashed, black lines) when $\gamma$ is high enough that condition~\eqref{condition1} is satisfied (i.e. condition~\eqref{condition} does not hold). Here, $\alpha_T=0.5$ and all the other parameters are as in Table~\ref{Tab:Tcr} with $\gamma=3.5$, $\eta=1.8$ and $\theta_C=\theta_T=1.8$. Panels \textbf{b.}-\textbf{d.} display similar plots for sufficiently large, intermediate and sufficiently small values of $\gamma$ that satisfy condition~\eqref{condition} -- i.e. $\gamma=2$ (panel \textbf{b.}), $\gamma=0.3$ (panel \textbf{c.}) and $\gamma=0.12$ (panel \textbf{d.}). All the other parameters are as in Table~\ref{Tab:Tcr} with $\eta=1.8$ and $\theta_C=\theta_T=1.8$. Initial conditions~\eqref{initcondAgentbased} and~\eqref{initcondPDE} with $a=1$ and $A=5$ were used to carry out numerical simulations. Analogous results were obtained when using different values of the parameter $A$ (results not shown). The results from the individual-based model correspond to the average over two realisations and the related variance is displayed by the coloured areas surrounding the curves.}
\label{fig:new1bis}
\end{figure}

	    			 \end{appendices}
   			\nocite{*}
\bibliography{bibliography}  	   			
\end{document}